\documentclass[12pt]{amsart} \usepackage{latexsym, amssymb}
\usepackage[T1]{fontenc}

\newtheoremstyle{slthm}% name
  {9pt}%      Space above, empty = `usual value'
  {5pt}%      Space below
  {\slshape}% Body font
  {}%         Indent amount (empty = no indent, \parindent = para indent)
  {\bfseries}% Thm head font
  {.}%        Punctuation after thm head
  {.5em}%     Space after thm head: " " = normal interword space;
  {\thmname{#1}\thmnumber{ #2}\thmnote{ (#3)}}% Thm head spec

\newtheoremstyle{prcl}% name
  {9pt}%      Space above, empty = `usual value'
  {5pt}%      Space below
  {\slshape}% Body font
  {}%         Indent amount (empty = no indent, \parindent = para indent)
  {\bfseries}% Thm head font
  {.}%        Punctuation after thm head
  {.5em}%     Space after thm head: " " = normal interword space;
  {\thmname{#3}\thmnumber{ #2}}% Thm head spec

\theoremstyle{slthm}
\newtheorem{thm}{Theorem}[section]
\newtheorem{lemma}[thm]{Lemma}
\newtheorem{prop}[thm]{Proposition}
\newtheorem{cor}[thm]{Corollary}

\theoremstyle{definition}

\newtheorem{df}[thm]{Definition}

\newtheorem{nrmk}[thm]{Remark}
\newtheorem{nrmks}[thm]{Remarks}
\newtheorem{expl}[thm]{Example}
\newtheorem{expls}[thm]{Examples}

\newtheorem*{term}{Terminology}

\theoremstyle{remark}
\newtheorem*{rmk}{Remark}

\theoremstyle{prcl}
\newtheorem*{prclaim}{Proclaim}

\newenvironment{renumerate}
        {\renewcommand{\labelenumi}{\rm (\roman{enumi})}
         \begin{enumerate}}
        {\end{enumerate}}

\newenvironment{lenumerate}[2]
        {\renewcommand{\labelenumi}{\rm (#1\theenumi)}
         \begin{enumerate}{\setcounter{enumi}{#2}}}
        {\end{enumerate}}

\newcounter{flexnummark}

\DeclareMathOperator{\length}{length}
\DeclareMathOperator{\cl}{cl}
\DeclareMathOperator{\fr}{fr}
\DeclareMathOperator{\bd}{bd}
\DeclareMathOperator{\ir}{int}

\DeclareMathOperator\gr{gr}

\DeclareMathOperator\h{h}

\DeclareMathOperator\fix{Fix}
\DeclareMathOperator\Lim{Lim}
\DeclareMathOperator\fixbd{Bd}
\DeclareMathOperator\urank{U^{\mbox{{\tiny\th}}}}

\newcommand{\rest}[1]{|_{#1}}

\newcommand{\into}{\longrightarrow}

\renewcommand{\bar}{\overline}

\newcommand{\acl}{\mbox{acl}}
\newcommand{\dcl}{\mbox{dcl}}

\newcommand{\set}[1]{\left\{#1\right\}}

\newcommand{\NN}{\mathbb{N}}

\newcommand{\RR}{\mathbb{R}}

\newcommand{\curly}[1]{\mathcal{#1}}

\newcommand{\C}{\curly{C}}

\newcommand{\F}{\curly{F}}
\newcommand{\G}{\curly{G}}

\newcommand{\M}{\curly{M}}

\renewcommand{\O}{\curly{O}}

\newcommand{\R}{\curly{R}}

\renewcommand{\b}{\mathfrak{b}}
\newcommand{\f}{\mathfrak{f}}
\newcommand{\g}{\mathfrak{g}}
\renewcommand{\h}{\mathfrak{h}}

\newcommand{\la}{\curly{L}}

\newcommand{\Creg}{\C_{\textrm{reg}}}
\newcommand{\Copen}{\C_{\textrm{open}}}
\newcommand{\Ctan}{\C_{\textrm{tan}}}
\newcommand{\Ctrans}{\C_{\textrm{trans}}}
\newcommand{\Csingle}{\C_{\textrm{single}}}
\newcommand{\phiopen}{\Phi_{\textrm{open}}}
\newcommand{\phitan}{\Phi_{\textrm{tan}}}
\DeclareMathOperator{\phitrans}{\Phi_{\text{trans}}}
\newcommand{\phisingle}{\Phi_{\textrm{single}}}

\numberwithin{equation}{section}

\allowdisplaybreaks[2]

\title {An ordered structure of rank two related to Dulac's
  Problem}

\author {A. Dolich and P. Speissegger}

\address {Department of Mathematics, Chicago State University, 
Chicago IL, USA}

\email{adolich@scu.edu}

\address {Department of Mathematics \& Statistics, McMaster
University, Hamilton ON, Canada}

\email {speisseg@math.mcmaster.ca}

\date {\today}



\thanks {Supported in parts by NSERC and NSF}

\begin{document}

\begin{abstract}
  For a vector field $\xi$ on $\RR^2$ we construct, under certain
  assumptions on $\xi$, an ordered model-theoretic structure
  associated to the flow of $\xi$.  We do this in such a way that the
  set of all limit cycles of $\xi$ is represented by a definable set.
  This allows us to give two restatements of Dulac's Problem for
  $\xi$---that is, the question whether $\xi$ has finitely many limit
  cycles---in model-theoretic terms, one involving the recently
  developed notion of $\urank$-rank and the other involving the notion
  of o-minimality.
\end{abstract}

\maketitle
\markboth{A. DOLICH AND P. SPEISSEGGER}{AN ORDERED STRUCTURE OF RANK TWO}

\section*{Introduction}

Let $\xi = a_1 \frac{\partial}{\partial x} + a_2
\frac{\partial}{\partial y}$ be a vector field on $\RR^2$ of class
$C^1$, and let
\begin{equation*}
  S(\xi) := \set{(x,y) \in \RR^2:\ a_1(x,y) = a_2(x,y) = 0}
\end{equation*}
be the set of singularities of $\xi$.  By the existence and uniqueness
theorems for ordinary differential equations (see Camacho and Lins
Neto \cite[p. 28]{cam-lin:foliations} for details), $\xi$ induces a
$C^1$-foliation $\F^\xi$ on $\RR^2 \setminus S(\xi)$ of dimension $1$.
Abusing terminology, we simply call a leaf of this foliation a
\textbf{leaf of $\xi$}.  A \textbf{cycle} of $\xi$ is a compact leaf
of $\xi$; a \textbf{limit cycle} of $\xi$ is a cycle $L$ of $\xi$ for
which there exists a non-compact leaf $L'$ of $\xi$ such that $L$ is
contained in the closure of $L'$.  

Dulac's Problem is the following statement: ``if $\xi$ is polynomial,
then $\xi$ has finitely many limit cycles''.  It is a weakening of the
second part of Hilbert's 16th problem, which states that ``there is a
function $H: \NN \into \NN$ such that for all $d \in \NN$, if $\xi$ is
polynomial of degree $d$ then $\xi$ has at most $H(d)$ limit cycles''.
Both problems have an interesting history, and while Dulac's problem
was independently settled in the 1990s by Ecalle \cite{eca:dulac} and
Ilyashenko \cite{ily:dulac}, Hilbert's 16th problem remains open; see
\cite{ily:dulac} for more details.

In this paper, we attempt to reformulate Dulac's Problem in
model-theoretic terms.  Our motivation to do so is twofold: we want to
\begin{renumerate}
\item find a model-theoretic structure naturally associated to $\xi$
  in which the flow of $\xi$ and the set of limit cycles of $\xi$ are
  represented by definable sets;
\item know to what extent the geometry of such a structure is
  determined by Dulac's Problem.
\end{renumerate}

Our starting point for (i) is motivated by the piecewise triviality of
Rolle foliations associated to analytic $1$-forms as described by
Chazal \cite{cha:rolle}.  Let $U \subseteq \RR^2$ be open; a leaf $L$
of $\xi\rest{U}$ is a \textbf{Rolle leaf} of $\xi\rest{U}$ if for
every $C^1$-curve $\delta:[0,1] \into U$ with $\delta(0) \in L$ and
$\delta(1) \in L$, there is a $t \in [0,1]$ such that $\delta'(t)$ is
tangent to $\xi(\delta(t))$.  Based on Khovanskii theory
\cite{kho:fewnomials} over an o-minimal expansion of the real field
\cite{spe:pfaffian}, we establish (Proposition \ref{constant_ranks}
and Theorem \ref{piecewise_trivial}):

\begin{prclaim}[Theorem A]
  Assume that $\xi$ is definable in an o-minimal expansion of the real
  field.  Then there is a cell decomposition $\C$ of $\,\RR^2$
  compatible with $S(\xi)$ such that, with $\Creg := \{C \in \C:\ C
  \cap S(\xi) = \emptyset\}$,
  \begin{enumerate}
  \item every $1$-dimensional $C \in \Creg$ is either transverse to
    $\xi$ or tangent to $\xi$;
  \item for every open $C \in \Creg$, every leaf of $\xi \rest{C}$ is a
    Rolle leaf of $\xi\rest{C}$;
  \item for every open $C \in \Creg$, the flow of $\xi$ in $C$ is
    represented by a lexicographic ordering of $C$.
  \end{enumerate}
\end{prclaim}

Part (3) of this theorem needs some explanation, as it represents our
understanding of the ``triviality'' of the flow of $\xi$ in $C$.
Given an open $C \in \Creg$, it follows from part (2) that the
direction of $\xi$ induces a linear ordering $<_\Gamma$ on every leaf
$L$ of $\xi\rest{C}$.  We can furthermore define a relation on the set
$\la(C)$ of all leaves of $\xi\rest{C}$ as follows: given a leaf $L$
of $\xi\rest{C}$, the fact that $L$ is a Rolle leaf of $\xi\rest{C}$
implies (see Remark \ref{hirsch} below) that $L$ separates $C
\setminus L$ into two connected components $U_{L,1}$ and $U_{L,2}$
such that the vector $\xi^\perp(z) := (a_2(z), -a_1(z))$ points into
$U_{L,2}$ for all $z \in L$.  Thus, for a leaf $L'$ of $\xi\rest{C}$
different from $L$, we define $L \ll_C L'$ if $L' \subseteq U_{L,2}$
and $L' \ll_C L$ if $L' \in U_{L,1}$.  In general, though, the
relation $\ll_C$ does not always define an ordering, even if every
leaf of $\xi\rest{C}$ is Rolle; see Example \ref{non-order} below.

Part (3) now means that the cell decomposition $\C$ may be chosen in
such a way that for every open $C \in \Creg$, the ordering $\ll_C$ on
$\la(C)$ is a linear ordering.  (See Example \ref{radial} for such a
decomposition in the situation of Example \ref{non-order}.)  This
leads to lexicographic orderings as follows: given $C \in \Creg$ and
$z \in C$, we denote by $L_z$ the leaf of $\xi\rest{C}$ containing
$z$.  If $C \in \Creg$ is open, we define a linear ordering $<_C$ on
$C$ by $x <_C y$ if and only if either $L_x \ll_C L_y$, or $L_x = L_y$
and $x <_{L_x} y$.  Letting $E_C$ be a set of representatives of
$\la(C)$, it is not hard to see that the structures $(C,<_C,E_C)$ and
$(\RR^2,<_{\text{lex}}, \{y=0\})$ are isomorphic, where
$<_{\text{lex}}$ is the usual lexicographic ordering of $\RR^2$.

To complete the picture, we also define an ordering $<_C$ on each
$1$-dimen\-sio\-nal $C \in \Creg$: if $C$ is tangent to $\xi$, we let
$<_C$ be the linear ordering induced on $C$ by the direction of $\xi$,
and if $C$ is transverse to $\xi$, we let $<_C$ be the linear ordering
induced on $C$ by the direction of $\xi^\perp$.  For each open $C \in
\Creg$, we also let $<_{E_C}$ be the restriction of $<_C$ to $E_C$.
Each of these orderings induces a topology on the corresponding set
that makes it homeomorphic to the real line.  Finally, for each
$1$-dimensional $C \in \Creg$ tangent to $\xi$, we fix an element $e_C
\in C$. 

In the situation of Theorem A, we reconnect the pieces of $\C$
according to the flow of $\xi$ as follows: let $B$ be the union of 
\begin{itemize}
\item all $1$-dimensional cells in $\Creg$ transverse to $\xi$,
\item the sets $E_C$ for all open cells $C \in \Creg$,
\item all $0$-dimensional cells in $\Creg$, and
\item the singletons $\{e_C\}$ for all $1$-dimensional $C \in \Creg$
  tangent to $\xi$.  
\end{itemize}
We define the \textbf{forward progression map} $\f:B \cup \{\infty\}
\into B \cup \{\infty\}$ by (roughly speaking) putting $\f(x)$ equal
to the next point in $B$ on the leaf of $\xi$ through $x$ if $x \ne
\infty$ and if such a point exists, and otherwise we put $\f(x) :=
\infty$.  In this situation, a point $x \in B$ belongs to a cycle of
$\xi$ if and only if there is a nonzero $n \in \NN$ such that $\f^n(x)
= x$, where $\f^n$ denotes the $n$-th iterate of $\f$.

In fact, only finitely many iterates of $\f$ are necessary to capture
all cycles of $\xi$ (Proposition \ref{fixed_points}): since a cycle of
$\xi$ is a Jordan curve in $\RR^2$, it is a Rolle leaf of $\xi$ and
therefore intersects each $C \in \C$ of dimension at most $1$ in at
most one connected component.  Hence there is an $N \in \NN$ such that
for all $x \in B$, $x$ belongs to a cycle of $\xi$ if and only if
$\f^N(x) = x$.

To see how we can use this to detect limit cycles of certain $\xi$, we
first define a cycle $L$ of $\xi$ to be a \textbf{boundary cycle}, if
for every $x \in L$ and every neighborhood $V$ of $x$, the set $V$
intersects some non-compact leaf of $\xi$.  Boundary cycles and limit
cycles are the same if $\xi$ is real analytic, because of the
following theorem of Poincar\'e's \cite{poincare} (see also Perko
\cite[p. 217]{perko}): 

\begin{prclaim}[Fact]
  If $\xi$ is real analytic, then $\xi$ cannot have an infinite number
  of limit cycles that accumulate on a cycle of $\xi$.
\end{prclaim}

On the other hand, it follows from the previous paragraph that for
every $x \in B$, the point $x$ belongs to a boundary cycle of $\xi$ if
and only if $x$ is in the boundary (relative to $B$ considered with
the topology induced on it by the various orderings defined above) of
the set of all fixed points of $\f^N$.

Based on the observations mentioned in the preceding paragraphs (and a
few related observations), we associate to each decomposition $\C$ as
in Theorem A a \textbf{flow configuration} $\Phi_\xi = \Phi_\xi(\C)$
of $\xi$, intended to code how the cells in $\C$ are linked together
by the flow of $\xi$.  To each flow configuration $\Phi$, we associate
in turn a unique first-order language $\la(\Phi)$, in such a way that
the situation described in the preceding paragraphs naturally yields
an $\la(\Phi_\xi)$-structure $\M_\xi$ in which the lexicographic
orderings of Theorem A, the associated forward progression map $\f:B
\cup \{\infty\} \into B \cup \{\infty\}$ and the set of all $x \in B$
that belong to some boundary cycle of $\xi$ are definable.

If, in the situation of Theorem A, there is an open $C \in \Creg$,
then the induced structure on $C$ in $\M_\xi$ is not o-minimal
(because the structure $(C,<_C,E_C)$ described above is definable in
$\M_\xi$).  Thus, to answer (ii) we need to work with notions weaker
than o-minimality.  A weakening that includes lexicographic
orderings is provided by the \textit{rosy} theories introduced
by Onshuus \cite{alf2}.

To recall this rather technical definition, we fix a complete first
order theory $T$ and a sufficiently saturated model $\M$ of $T$, and
we work in $\M^{\rm{eq}}$.  (For standard model-theoretic terminology, we
refer the reader to Marker \cite{mark2}).  The definition of \th
-forking is much like that of forking in the stable or simple context:
A formula $\phi(x,a)$ \textbf{strongly divides} over a set $A$ if
tp$(a/A)$ is non-algebraic and the set $\{\phi(x,b): b \models$
tp$(a/A)\}$ is $k$-inconsistent for some $k \in \NN$.  The formula
$\phi(x,a)$ \textbf{\th -divides} over $A$ if for some tuple $c$,
$\phi(x,a)$ strongly divides over $A \cup \{c\}$.  The formula
$\phi(x,a)$ \textbf{\th -forks} over $A$ if $\phi(x,a)$ implies a
finite disjunction of formulas all of which \th -divide over $A$.  A
complete type $p(x)$ \textbf{\th -forks} over $A$ if there is some
formula $\phi(x)$ in $p(x)$ that \th -forks over $A$.

For a theory $T$ to be rosy means, roughly speaking, that in models of
$T$, \th -forking has many desirable properties, much like forking in
the stable or simple contexts.  For the formal definition we need only
focus on a single one of these: $T$ is \textbf{rosy} if for any
complete type $p(x)$ over a parameter set $B$, there exists $B_0
\subseteq B$ with $\| B_0 \| \leq \| T \|$ such that $p(x)$ does not
\th -fork over $B_0$.

The ``degree of rosiness'' of a theory is measured by the
$\urank$-rank, defined analogously to the $U$-rank in stable theories.
For an ordinal $\alpha$ and a complete type $p(x)$ with parameter set
$A$, we define $\urank(p) \geq \alpha$ by ordinal induction:
\begin{renumerate}
\item $\urank(p)\geq 0$ if $p$ is consistent;
\item if $\alpha$ is a limit ordinal, then $\urank(p)\geq \alpha$ if
  $\urank(p) \geq \beta$ for all $\beta < \alpha$;
\item $\urank(p) \geq \alpha+1$ if there is a complete type $q(x)$ so
  that $p \subseteq q$, $q$ \th -forks over $A$ and $\urank(q)\geq
  \alpha$.
\end{renumerate}
For an ordinal $\alpha$, we say that $\urank(p) = \alpha$ if
$\urank(p) \geq \alpha$ and $\urank(p) \not \geq \alpha+1$.  Finally,
$\urank(T)$ is defined to be the supremum of $\urank(p)$ for all
one-types $p$ with parameters over the empty set.  One of the
fundamental facts about rosy theories is that $T$ is rosy if
$\urank(T)$ is an ordinal \cite{alf2}.

For example, every o-minimal theory is rosy of $\urank$-rank one.  On
the other hand, the theory $T$ of the structure $(C,<_C,E_C)$ above has
$\urank$-rank at least two.  To see the latter, let $\M \models T$ be
$\aleph_1$-saturated and write $C_z := \{x \in C:\ z_1 <_C x <_C z_2$ for all $z_1, z_2 \in E_C$ such that $z_1 <_C z <_C z_2\}$.  Since $E_C^\M$ is a dense linear ordering without endpoints, there are infinitely many $a \in E_C^\M$ such that $a \notin \acl(\emptyset)$.  For any two such $a, b \in E_C^\M$, the fibers $C_a^\M$ and $C_b^\M$ are disjoint, infinite definable sets.  Hence $\urank(\M) \ge 2$.  

In this paper, we use the argument of the previous example to
establish lower bounds on $\urank$-rank for the theories we are
interested in.  For upper bounds, we need a special case of the
Coordinatization Theorem \cite[Theorem 2.2.2]{alf1}:

\begin{prclaim}[Fact]
  Assume that $T$ defines a dense linear ordering without endpoints,
  and let $\M \models T$ be saturated.  Let also $n \in \NN$ and
  assume that for all $a \in M$, there are $a_1, \dots, a_n \in M$
  such that $a = a_n$ and for each $i \in \{1, \dots, n\}$, the type
  of $(a_1, \dots, a_i)$ over $(a_1, \dots, a_{i-1})$ is implied in
  $T$ by the order type of $(a_1, \dots, a_i)$ over $(a_1, \dots,
  a_{i-1})$.  Then $\urank(T) \le n$.
\end{prclaim}

Note that our discussion above and the previous example imply that
$\urank(\M_\xi) \ge 2$.  The main result of this paper is the
following restatement of Dulac's problem:

\begin{prclaim}[Theorem B]
  Assume that $\xi$ is definable in an o-minimal expansion of the real
  field, and let $\M_\xi$ be the $\la(\Phi_\xi)$-structure associated
  to some flow configuration $\Phi_\xi$ of $\xi$.  Then
  \begin{enumerate}
  \item $\xi$ has finitely many boundary cycles if and only if
    $\,\urank(\M_\xi) = 2$;
  \item if $\xi$ is real analytic, then $\xi$ has finitely many limit
    cycles if and only if $\,\urank(\M_\xi) = 2$.
  \end{enumerate}
  \smallskip
\end{prclaim}

The proof of Theorem B is lengthy, but straightforward: we prove that
$\M_\xi$ admits quantifier elimination in a certain expanded language
(Theorem \ref{main}).  The main ingredient in this proof is a
reduction---modulo the theory of $\M_\xi$ in the expanded language,
roughly speaking---of general quantifier-free formulas to certain
quantifier-free order formulas, which allows us to deduce the
quantifier elimination for $\M_\xi$ from quantifier elimination of the
theory of $(\RR^2, <_{\text{lex}}, \{y=0\}, \pi)$, where $\pi:\RR^2
\into \{y=0\}$ is the canonical projection on the $x$-axis.  Under the
assumption of having only finitely many boundary cycles, the new
predicates of the expanded language are easily seen to define subsets
of the various cells obtained by Theorem A that are finite unions of
points and intervals.  Sufficiency in Theorem B then follows from the
above Fact; necessity follows by general $\urank$-rank arguments.

As a corollary of Theorem B, Ecalle's and Ilyashenko's solutions of
Dulac's Problem imply the following:

\begin{prclaim}[Corollary]
  Assume $\xi$ is polynomial, and let $\M_\xi$ be the
  $\la(\Phi_\xi)$-structure associated to some flow configuration
  $\Phi_\xi$ of $\xi$.  Then $\urank(\M_\xi) = 2$. \qed
\end{prclaim}

It remains an open question whether, in the situation of the corollary,
the structures are definable in some o-minimal expansion of the real
line.  An answer to this question, however, seems to go far beyond our
current knowledge surrounding Dulac's Problem.

Finally, our proof of Theorem B gives rise to a second restatement of
Dulac's problem that does not involve $\urank$-rank: let $G$ be the
union of all $1$-dimensional $C \in \Creg$ that are transverse to
$\xi$, all $0$-dimensional $C \in \Creg$ and $\{\infty\}$.  Let
$\G_\xi$ be the expansion of $G$ by all corresponding orderings $<_C$
and by the map $\f^2 \rest{G}$.  (Note that $\f^2\rest{G}$ maps $G$
into $G$.)  We may view $\G_\xi$ as a graph whose vertices are the
elements of $G$ and whose edges are defined by $\f^2$.

\begin{prclaim}[Theorem C]
  Assume that $\xi$ is definable in an o-minimal expansion of the real
  field, and let $\G_\xi$ be as above.  Then 
  \begin{enumerate}
  \item $\xi$ has finitely many boundary cycles if and only if the
    structure induced by $\G_\xi$ on each $1$-dimensional $C \subseteq
    G$ is o-minimal;
  \item if $\xi$ is real analytic, then $\xi$ has finitely many limit
    cycles if and only if the structure induced by $\G_\xi$ on each
    $1$-dimensional $C \subseteq G$ is o-minimal.
  \end{enumerate} 
\end{prclaim}

Our paper is organized as follows: in Sections
\ref{rolle}--\ref{piecewise_trivial_section}, we establish Theorem A:
in Section \ref{rolle}, we combine basic o-minimal calculus with
Khovanskii's Lemma to obtain a cell decomposition satisfying (1) and
(2) of Theorem A.  To refine this decomposition so that (3) holds, we
need to study what sets we obtain as Hausdorff limits of a sequence of
leaves of $\xi \rest{C}$ (Proposition \ref{H-limit}).  The refinement
is then given in Section \ref{piecewise_trivial_section}, where (3) is
established as Theorem \ref{piecewise_trivial}.  In Sections
\ref{foliation_orderings} and \ref{progression_map}, we define the
relevant orderings and progression maps associated to $\xi$ as
mentioned earlier.  Inspired by the latter, we then introduce the
notion of a flow configuration and the associated first-order language
in Section \ref{flowconfig}, where we also give an axiomatization of
the crucial properties satisfied by the models $\M_\xi$ above.  Some
basic facts about the iterates of the forward progression map are
deduced from these axioms in Section \ref{towards}.  In Section
\ref{dulac_flow}, we extend our axioms to reflect the additional
assumption that there are only finitely many boundary cycles, and we
introduce additional predicates for certain definable sets related to
the sets of fixed points of the iterates of the forward progression
map.  The quantifier elimination result is then given in Section
\ref{QE}, and we prove Theorems C and B in Section \ref{consequences}.
We finish with a few questions and remarks in Section \ref{final}.
\smallskip

\noindent\textbf{Acknowledgements.} We thank Lou van den Dries and
Chris Miller for their suggestions and comments on the earlier
versions of this paper. \smallskip

\noindent\textbf{Global conventions.}  We fix an o-minimal expansion
$\R$ of the real field; ``definable'' means ``definable in $\R$ with
parameters''.  

For $1 \le m \le n$, we denote by $\Pi_m:\RR^n \into \RR^m$ the
projection on the first $m$ coordinates.  

Given $(x,y) \in \RR^2$, we put $(x,y)^\perp := (y,-x)$.

For a subset $A \subseteq \RR^n$, we let $\cl(A)$, $\ir(A)$, $\bd(A)
:= \cl(A) \setminus \ir(A)$ and $\fr(A) := \cl(A) \setminus A$ denote
the topological closure, interior, boundary and frontier,
respectively.

For $n \in \NN$, we define the analytic diffeomorphism $\phi_n:\RR^n
\into (-1,1)^n$ by $\phi_n(x_1, \dots, x_n):=
\left(x_1/\sqrt{1+x_1^2}, \dots, x_n/\sqrt{1+x_n^2}\right).$  Given $X
\subseteq \RR^n$, we write $X^*:= \phi_n(X)$, and given a vector field
$\eta$ on $\RR^n$ of class $C^1$, we write $\eta^*$ for the
push-forward $(\phi_n)_* \,\eta$ of $\eta$ to $(-1,1)^n$.

\section{Rolle decomposition}  \label{rolle}

Let $U \subseteq \RR^2$ be open and $p \ge 1$ be an integer.  Let $\xi
= a_1 \frac{\partial}{\partial x} + a_2 \frac{\partial}{\partial y}$
be a definable vector field on $U$ of class $C^p$ (that is, the
functions $a_1,a_2:U \into \RR$ are definable and of class $C^p$), and
let
\begin{equation*}
  S(\xi):= \set{z \in U:\ a_1(z) = a_2(z) = 0}
\end{equation*} 
be the set of singularities of $\xi$.  By the existence and uniqueness
theorems for ordinary differential equations
\cite[p. 28]{cam-lin:foliations}, $\xi$ induces a $C^p$-foliation
$\F^\xi$ on $U \setminus S(\xi)$ of dimension $1$.  Abusing
terminology, we simply call a leaf of this foliation a \textbf{leaf of
  $\xi$}.

\begin{rmk}
  Put $\omega := a_2 dx - a_1 dy$; then $S(\xi)$ is the set of
  singularities of $\omega$, and the foliation $\F^\xi$ is exactly the
  foliation on $U \setminus S(\xi)$ defined by the equation $\omega =
  0$.  Below, we will use this observation (mainly in connection with
  some citations) without further mention.
\end{rmk}

\begin{df}  
  \label{piecewise_transverse}
  Let $\gamma:I \into U$ of class $C^p$, where $I \subseteq \RR$ is an
  interval.  We call $\gamma$ a \textbf{$C^p$-curve in $U$} and
  usually write $\Gamma:= \gamma(I)$.  If $t \in I$ is such that
  $\xi^\perp(\gamma(t))\cdot\gamma'(t) \ne 0$, we say that $\gamma$ is
  \textbf{transverse to} $\xi$ \textbf{at $t$}; otherwise, $\gamma$ is
  \textbf{tangent to $\xi$ at $t$}.  The curve $\gamma$ is
  \textbf{transverse (tangent) to $\xi$} if $\gamma$ is transverse
  (tangent) to $\xi$ at every $t \in I$.

  A leaf $L$ of $\xi$ is a \textbf{Rolle leaf of $\xi$} if for every
  $C^1$-curve $\gamma:[0,1] \into U$ with $\gamma(0) \in L$ and
  $\gamma(1) \in L$, there is a $t \in [0,1]$ such that
  $\xi^\perp(\gamma(t)) \cdot \gamma'(t) = 0$.

  A \textbf{cycle} of $\xi$ is a compact leaf of $\xi$.  A cycle
  $L$ of $\xi$ is a \textbf{limit cycle} of $\xi$ if there is a
  non-compact leaf $L'$ of $\xi$ such that $L \subseteq \cl(L')$.
  A cycle $L$ of $\xi$ is a \textbf{boundary cycle} of $\xi$ if
  for every open set $V \subseteq \RR^2$ with $V \cap L \ne
  \emptyset$, there is a non-compact leaf $L'$ of $\xi$ such that
  $V \cap L' \ne \emptyset$.
\end{df}

\begin{nrmk}
  \label{hirsch}
  Since $\xi$ is integrable in $U \setminus S(\xi)$, every Rolle leaf
  $L$ of $\xi$ is an embedded submanifold of $U \setminus S(\xi)$ that
  is closed in $U \setminus S(\xi)$.  In particular, by Theorem 4.6
  and Lemma 4.4 of Chapter 4 in \cite{hir:difftop}, if $U \setminus
  S(\xi)$ is simply connected, then $U \setminus (S(\xi) \cup L)$ has
  exactly two connected components such that $L$ is equal to the
  boundary in $U \setminus S(\xi)$ of each of these components.
\end{nrmk}

\begin{lemma}[Khovanskii \cite{kho:fewnomials}]  
  \label{separating}
  \begin{enumerate}
  \item Assume that $U \setminus S(\xi)$ is simply connected, and
    let $L \subseteq U\setminus S(\xi)$ be an embedded leaf of
    $\xi$ that is closed in $U \setminus S(\xi)$.  Then $L$ is a
    Rolle leaf of $\xi$ in $U$.
  \item Let $L$ be a cycle of $\xi$.  Then $L$ is a Rolle leaf of
    $\xi$.
  \end{enumerate}
\end{lemma}

\begin{proof}[Sketch of proof]
  (1) Arguing as in the preceding remark, the set $U \setminus S(\xi)$
  has exactly two connected components $U_1$ and $U_2$, such that
  $\bd(U_i) \cap (U \setminus S(\xi)) = L$ for $i=1,2$.  The argument
  of Example 1.3 in \cite{spe:pfaffian} now shows that $L$ is a Rolle
  leaf of $\xi$.

  (2) Since $L$ is compact, $L$ is an embedded and closed submanifold
  of $\RR^2$.  Now conclude as in part (1).
\end{proof}

\begin{df}  
  \label{rolle_form}
  We call $\xi$ \textbf{Rolle} if $S(\xi) = \emptyset$, $\xi$
  is of class $C^1$ and every leaf of $\xi$ is a Rolle leaf of
  $\xi$.
\end{df}

We now let $\C$ be a $C^p$-cell decomposition of $\RR^2$ compatible
with $U$ and $S(\xi)$, and we put $\C_U:= \set{C \in \C:\ C
  \subseteq U}$.  Refining $\C$, we may assume that $\xi\rest{C}$
is of class $C^p$ for every $C \in \C_U$, and that every $C \in \C_U$
of dimension $1$ is either tangent or transverse to $\xi$.
Refining $\C$ again, we also assume that
\begin{itemize}
\item[(I)] $a_1$ and $a_2$ have constant sign on every $C \in \C_U$.
\end{itemize}
Such a decomposition $\C$ is called a \textbf{Rolle decomposition for
  $\xi$}, because of the following:

\begin{prop}  
  \label{constant_ranks}
  Let $C \in \C_U$ be open such that $C \cap S(\xi) = \emptyset$.
  Then $\xi\rest{C}$ is Rolle.  Moreover, if both $a_1$ and $a_2$
  have nonzero constant sign on $C$, then either every leaf of
  $\xi\rest{C}$ is the graph of a strictly increasing $C^p$
  function $f:I \into \RR$, or every leaf of $\xi\rest C$ is the
  graph of a strictly decreasing $C^p$-function $f:I \into \RR$, where
  $I \subseteq \RR$ is an open interval depending on $f$.
\end{prop}

\begin{proof}
  If $a_1\rest{C} = 0$ or $a_2\rest{C} = 0$, the conclusion is
  obvious.  So we assume that $a_1\rest{C}$ and $a_2\rest{C}$ have
  constant positive sign, say; the remaining three cases are handled
  similarly.  Let $L$ be a leaf of $\xi\rest{C}$; we claim that $L$ is
  the graph of a strictly increasing $C^p$-function $f:I \into \RR$,
  where $I:= \Pi_1(L)$.

  To see this, assume first that there are $x,y_1,y_2 \in \RR$ such
  that $(x,y_i) \in L$ for $i=1,2$ and $y_1 \ne y_2$.  Since
  $\xi\rest{C}$ is of class $C^p$, the leaf $L$ is a $C^p$-curve, so
  by Rolle's Theorem, there is an $a \in L$ such that $L$ is tangent
  at $a$ to $\partial / \partial y$.  But this means that $a_1(a) =
  0$, a contradiction.  Thus, $L$ is the graph of a strictly
  increasing $C^p$-function $f:I \into \RR$.

  It follows from the claim that $L$ is an embedded submanifold of $C$
  and, since $C \cap S(\xi) = \emptyset$, that $L$ is a closed
  subset of $C$.  Thus by Lemma \ref{separating}(1), $L$ is a Rolle
  leaf of $\xi\rest C$.
\end{proof}

\section{Rolle foliations and Hausdorff limits of Rolle leaves}
\label{hausdorff}

We continue working with $\xi$ as in Section \ref{rolle}, and we fix a
Rolle decomposition $\C$ for $\xi$.  We fix an open $C \in \C_U$ such
that $C \cap S(\xi) = \emptyset$.

To simplify notation, we write $\xi$ in place of $\xi\rest{C}$
throughout this section.

Let $L$ be a leaf of $\xi$.  Since $L$ is a Rolle leaf of $\xi$, $C
\setminus L$ has two connected components $U_{L,1}$ and $U_{L,2}$, and
$L$ is the boundary of $U_{L,i}$ in $C$ for $i=1,2$.  Since $\xi^\perp(z)
\ne (0,0)$ for all $z \in C$ and $L$ is connected, there is an $i \in
\{1,2\}$ such that $\xi^\perp(z)$ points inside $U_{L,i}$ for all $z \in L$;
reindexing if necessary, we may assume that $\xi^\perp(z)$ points inside
$U_{L,2}$ for every leaf $L$ of $\xi$.

\begin{df}  
  \label{foliation}
  For a point $z \in C$, we let $L^\xi_z$ be the unique leaf of
  $\xi$ such that $z \in L^\xi_z$.  For any subset $X \subseteq
  C$, we define
  \begin{equation*}
    F^\xi(X):= \bigcup_{z \in X} L^\xi_z,
  \end{equation*}
  called the \textbf{$\xi$-saturation of $X$}, and we put
  \begin{equation*}
    \la^\xi(X):= \set{L^\xi_z:\ z \in X}.
  \end{equation*} 
  For $X \subseteq C$, we define a relation $\ll_X^\xi$ on the set
  $\la^\xi(X)$ as follows: $L \ll_X^\xi M$ if and only if $L \subseteq
  U_{M,1}$ (if and only if $M \subseteq U_{L,2}$).

  Whenever $\xi$ is clear from context, we omit ``$\xi$'' in the
  definitions and notations above.
\end{df}

Note that in general the relation $\ll_C$ may not
define an order relation on $\la(C)$:

\begin{expl}
  \label{non-order}
  Let $\zeta := -y\frac{\partial}{\partial x} + x
  \frac{\partial}{\partial y}$, and let $g:\RR^2 \into \RR$ be defined
  by $g(x,y) := (y - (x-2))^2$.  Then $g\zeta$ is a real analytic
  vector field on $\RR^2$ and $S(g\zeta) = \{0\} \cup \{(x,y):\ y =
  x-1\}$.  Let also $C$ be the cell $(\alpha,\beta)$, where $\alpha,
  \beta:(0,1) \into \RR$ are defined by $\alpha(x) := x-2$ and
  $\beta(x) := x-1$.

  Then $C \cap S(g\zeta) = \emptyset$, and since every leaf of
  $\zeta$ is a Rolle leaf of $\zeta$, the vector field $g\zeta
  \rest{C}$ is Rolle.  However, $\ll_C^{g\zeta}$ is not an ordering of
  $\la(C)$: pick a leaf $L$ of $\xi$ (that is, a circle with center
  $(0,0)$) such that $L \cap \gr(\alpha)$ contains two points.  Then
  $L \cap C$ consists of two distinct leaves $L_1$ and $L_2$ of
  $g\zeta \rest{C}$.  Since $\zeta^\perp(z)$ points outside the circle
  $L$ for every $z \in L$, we get $L_1 \subseteq U_{L_2,1}$ and $L_2
  \subseteq U_{L_1,1}$, that is, $L_1 \ll_C^{g\zeta} L_2$ and $L_2
  \ll_C^{g\zeta} L_1$.
\end{expl}

However, for certain $X$ the relation $\ll_X$ is a linear ordering of
$\la(X)$, as discussed in the following lemma.  For a curve $\gamma:I
\into C$, we write
\begin{equation*}
  L(t) := L_{\gamma(t)} \qquad \text{for all } t \in I;
\end{equation*}
in this situation, we have $F(\Gamma) = \bigcup_{t \in I} L(t)$.

\begin{lemma}  
  \label{order_along_curve}
  Let $\gamma:I \into C$ be a $C^p$-curve transverse to $\xi$,
  where $I \subseteq \RR$ is an interval.
  \begin{enumerate}
  \item If $I$ is open, then $F(\Gamma)$ is open.
  \item The relation $\ll_\Gamma$ is a linear ordering of
    $\la(\Gamma)$, and the map $t \mapsto L(t):I \into \la(\Gamma)$ is
    order-preserving if $\xi^\perp(\gamma(t)) \cdot \gamma'(t) > 0$
    for all $t \in I$ and order-reversing if $\xi^\perp(\gamma(t))
    \cdot \gamma'(t) < 0$ for all $t \in I$.
  \end{enumerate}
\end{lemma}

\begin{proof}
  (1) Assume that $I$ is open, and let $t \in I$.  Because $\xi$ is
  $C^p$ and nonsingular and $\gamma$ is transverse to $\xi$, by a
  variant of Picard's Theorem (see Theorem 8-2 of
  \cite{aus-mac:differentiable}), there is an open set $B_t \subseteq C$
  containing $\gamma(t)$ such that $B_t \subseteq F(\Gamma)$.  Put
  $B:= \bigcup_{t \in I} B_t$; then $\Gamma \subseteq B \subseteq
  F(\Gamma)$, so $F(\Gamma) = F(B)$.  Since $B$ is open, it follows
  from Theorem III.1 in \cite{cam-lin:foliations} that $F(\Gamma)$ is
  open.

  (2) Since $\gamma$ is transverse to $\xi$ and each $L(t)$ is
  Rolle, the map $t \mapsto L(t):I \into \la(\Gamma)$ is injective.
  It therefore suffices to show that either
  \begin{equation*}
    s < t \quad\Leftrightarrow\quad L(s) \ll_\Gamma L(t) \qquad \text{for
      all } s,t \in I,
  \end{equation*}
  or
  \begin{equation*}
    s < t \quad\Leftrightarrow\quad L(t) \ll_\Gamma L(s) \qquad \text{for
      all } s,t \in I.
  \end{equation*}
  Since $\gamma$ is transverse to $\xi$, the continuous map $t \mapsto
  \xi^\perp(\gamma(t)) \cdot \gamma'(t):I \into \RR$ has constant
  positive or negative sign.  Assume it has constant positive sign;
  the case of constant negative sign is handled similarly.  Then for
  every $t \in I$, the set
  \begin{equation*}
    \Gamma_{<t}:= \set{\gamma(s):\ s \in I, s<t}
  \end{equation*}
  is contained in $U_{L(t),1}$.  Hence $L(s) \subseteq U_{L(t),1}$ for
  all $s \in I$ with $s<t$, that is, $L(s) \ll_\Gamma L(t)$ for all $s
  \in I$ with $s < t$.  Similarly, $L(t) \ll_\Gamma L(s)$ for all $s
  \in I$ with $s > t$, and since $t \in I$ was arbitrary, the lemma
  follows.
\end{proof}

We assume for the rest of this section that $C$ is \textit{bounded}.
Let $\xi_C$ be the $1$-form on $C$ defined by
\begin{equation*} 
  \xi_C:= \frac {\xi\rest{C}}{\|\xi\rest{C}\|}.
\end{equation*}
Then $\xi_C$ is a bounded, definable $C^p$-map on $C$, so by
o-minimality, there is a finite set $F_C \subseteq \fr(C)$ such that
$\xi_C$ extends continuously to $\cl(C) \setminus F_C$; we denote
this continuous extension by $\xi_C$ as well.

Let $c, d \in \RR$ and $\alpha,\beta:(c,d) \into \RR$ be definable and
$C^p$ such that $C = (\alpha,\beta)$.  Because $C$
is bounded, the limits $\alpha(c):= \lim_{x \to c} \alpha(x)$,
$\alpha(d):= \lim_{x \to d} \alpha(x)$, $\beta(c):= \lim_{x \to c}
\beta(x)$ and $\beta(d):= \lim_{x \to d} \beta(x)$ exist in $\RR$.
The points of the set
\begin{equation*}
  V_C:= \set{(c,\alpha(c)), (d,\alpha(d)), (c,\beta(c)),
    (d,\beta(d))}
\end{equation*}
are called the \textbf{corners} of $C$.

\begin{expl}
  \label{non-order2}
  In Example \ref{non-order}, we have $F_C \subseteq V_C$ and both $g\zeta
  \cdot (\partial/\partial x)$ and $g\zeta \cdot (\partial/\partial
  y)$ have constant nonzero sign.  The next proposition shows that
  under the latter assumptions, the situation of Example
  \ref{non-order} is as bad as it gets.
\end{expl}

\begin{prop}  
  \label{H-limit}
  Suppose that $F_C \subseteq V_C$, $a_1\rest{C} \ne 0$ and
  $a_2\rest{C} \ne 0$.  Let $\gamma:[0,1] \into C$ be a $C^p$-curve
  transverse to $\xi$, and let $t_i \in (0,1)$ be such that $t_0 <
  t_1 < t_2 < \cdots$ and $t_i \to 1$.  Then the sequence
  $\big(\cl(L(t_i))\big)$ converges in the Hausdorff metric to a
  compact set $K:= \lim \cl(L(t_i)) \subseteq \cl(C)$, such that
  \begin{renumerate}
  \item $\Pi_1(K) = [a,b]$ with $c \le a < b \le d$;
  \item each component of $K \cap C$ is a leaf of $\xi$;
  \item $K \cap \Pi_1^{-1}(a,b) = \gr(f)$ for some continuous function
    $f:(a,b) \into \RR$.
  \end{renumerate}
\end{prop}

\begin{proof}
  By Proposition \ref{constant_ranks}, we may assume that for every $t
  \in [0,1]$, the leaf $L(t)$ is the graph of a strictly increasing
  $C^p$-function $f_t:(a(t),b(t)) \into \RR$ (the other cases are
  handled similarly).  Since $C$ is bounded, the limits $f_t(a(t)):=
  \lim_{x \to a(t)} f_t(x)$ and $f_t(b(t)):= \lim_{x \to b(t)} f_t(x)$
  exist, and we also denote by $f_t:[a(t),b(t)] \into \RR$ the
  corresponding continuous extension of $f_t$.  Then $\cl(L(t)) =
  \gr(f_t)$.  By Lemma \ref{order_along_curve}, we may also assume
  that the map $t \mapsto L(t):[0,1] \into \la(\Gamma)$ is
  order-preserving (again, the other case is handled similarly).
  Finally, since each $f_{t}$ is strictly increasing and the map $t
  \mapsto L(t):[0,1] \into \la(\Gamma)$ is order-preserving, it
  follows that $f_{s}(x) > f_{t}(x)$ for all $s,t \in [0,1]$ such that
  $s<t$ and $x \in (a(s),b(s)) \cap (a(t), b(t))$.

  Since each $\cl(L(t_i))$ is connected, the set $K$ is connected, so
  $\Pi_1(K)$ is an interval $[a,b]$, which proves (i).  It follows in
  particular that for every $x \in (a,b)$, there is an open interval
  $I_x \subseteq (a,b)$ containing $x$ such that $I_x \subseteq
  (a(t_i),b(t_i))$ for all sufficiently large $i$.  Thus by our
  assumptions,
  \begin{itemize}
  \item[($\ast$)] for every $x \in (a,b)$ we have $f_{t_i}\rest{I_x} >
    f_{t_{i+1}}\rest{I_x}$ for sufficiently large $i$.
  \end{itemize}

  Next, we show that $K \cap C$ is an integral manifold of $\xi$.  Fix
  a point $(x,y) \in K \cap C$; it suffices to show that there is an
  open box $B \subseteq C$ containing $(x,y)$ such that $K \cap B$ is
  an integral manifold of $\xi$.  Let $B = I \times J$ be an open box
  containing $(x,y)$ such that $I \subseteq I_x$.  Since $a_1(x,y) \ne
  0$, we may also assume (after shrinking $B$) that there is an
  $\epsilon > 0$ such that $|a_1(x',y')| \ge \epsilon$ for all
  $(x',y') \in B$; in particular, there is an $M>0$ such that
  $f_{t_i}\rest{I}$ is $M$-Lipshitz for all sufficiently large $i$.
  Hence by ($\ast$), the function $f:I \into \RR$ defined by $f(x'):=
  \lim_{i \to \infty} f_{t_i}(x')$ is Lipshitz and satisfies $K \cap
  (I \times \RR) = K \cap B = \gr(f)$.  Finally, shrinking $B$ again
  if necessary, the fact that $\F^\xi$ is a foliation gives that $K
  \cap B$ is an integral manifold of $\xi$, as required.

  Since $K$ is compact and $K \cap C$ is an integral manifold of
  $\xi$, every component of $K \cap C$ is a leaf of $\xi$.  It
  also follows from the previous paragraph that $K \cap C$ is the
  graph of a continuous function $g:\Pi_1(K \cap C) \into \RR$, which
  proves (ii).

  Let now $x \in (a,b)$ be such that $x \notin \Pi_1(K \cap C)$.  Then
  $(x,\alpha(x))$ or $(x,\beta(x))$ belongs to $K$, because $(a,b)
  \subseteq \Pi_1(K)$; by ($\ast$) we have $(x,\beta(x)) \notin K$, so
  $(x,\alpha(x)) \in K$.  If $(\xi_C \cdot \frac{\partial}{\partial
    x})(x,\alpha(x)) \ne 0$, then by the same arguments as used for
  (ii), we conclude that there are open intervals $I,J \subseteq \RR$
  such that $(x,\alpha(x)) \in I \times J$ and $K \cap (I \times J)$
  is the graph of a continuous function defined on $I$.  Therefore,
  part (iii) is proved once we show that $(\xi_C \cdot
  \frac{\partial}{\partial x})(x,\alpha(x)) \ne 0$ for all $x \in
  (a,b) \setminus \Pi_1(K \cap C)$.

  Assume for a contradiction that there is an $x \in (a,b) \setminus
  \Pi_1(K \cap C)$ such that $(\xi_C \cdot \frac{\partial}{\partial
    x})(x,\alpha(x)) = 0$.  Let $M > |\alpha'(x)|$, and let $I,J
  \subseteq \RR$ be open intervals such that $I \subseteq I_x$ and
  $|a_2/a_1| > M$ on $B:= I \times J$.  Since $f_{t_i}(x) \to
  \alpha(x)$, it follows from the fundamental theorem of calculus for
  all sufficiently large $i$ that $f_{t_i}(x_i) = \alpha(x_i)$ for some
  $x_i \in I$, a contradiction.
\end{proof}

\section{Piecewise trivial decomposition}  
\label{piecewise_trivial_section}

We continue working with $\xi$ as in Section \ref{rolle}, and we adopt
the notations used there.  Note that $\xi^*$ (as defined at the end of
the introduction) is a definable vector field on $U^*$ of class $C^p$,
and that $\C$ is a Rolle decomposition of $\RR^2$ for $\xi$ if and
only if $\C^*:= \set{C^*:\ C \in \C}$ is a Rolle decomposition of
$(-1,1)^2$ for $\xi^*$.

Let $C \subseteq U$ be a bounded, open, definable $C^p$-cell such that
$\xi\rest{C}$ is Rolle.  To detect situations like the one described
in Example \ref{non-order}, we associate the following notations to
such a $C$: there are real numbers $c < d$ and definable $C^p$
functions $\alpha,\beta:(c,d) \into \RR$ such that $C =
(\alpha,\beta)$.  Given a $C^1$-function $\delta:(c,d) \into \RR$ such
that $\alpha(x) \le \delta(x) \le \beta(x)$ for all $x \in (c,d)$, we
define $\sigma_\delta:C \into \RR$ by
\begin{equation*}
  \sigma_\delta(x,y) := \xi^\perp(x,y) \cdot \begin{pmatrix} 1 \\
    \delta'(x) \end{pmatrix}.
\end{equation*}
Note that for each $x \in (c,d)$, there are by o-minimality a maximal
$\alpha_0^C(x) \in (\alpha(x),\beta(x)]$ and a minimal $\beta_0^C(x)
\in [\alpha(x),\beta(x))$ such that the function $\sigma_\alpha$ has
constant sign on $\{x\} \times (\alpha(x),\alpha_0^C(x))$ and the
function $\sigma_\beta$ has constant sign on $\{x\} \times
(\beta_0^C(x),\beta(x))$; we omit the superscript ``$C$'' whenever $C$
is clear from context.  Note that $\alpha_0,\beta_0:(c,d) \into \RR$
are definable.

\begin{df} 
  \label{pi-tri} 
  A $C^p$-cell decomposition of $\RR^2$ compatible with $U$, $\bd(U)$
  and $S(\xi)$ is called \textbf{almost piecewise trivial for
    $\xi$} if
  \begin{itemize}
  \item[(I)] every $C \in \C_U$ of dimension $1$ is either tangent or
    transverse to $\xi$;
  \item[(II)] the components of $\xi$ have constant sign on every $C
    \in \C_U$;
  \end{itemize}
  and for every open, bounded $C \in \C_U$ such that $C \cap S(\xi)
  = \emptyset$, the following hold:
  \begin{itemize}
  \item[(III)] $F_C \subseteq V_C$;
  \item[(IV)] the maps $\alpha_0,\beta_0:(c,d) \into \RR$ are
    continuous;
  \item[(V)] the map $\sigma_\alpha$ has constant sign on the cell
    $(\alpha,\alpha_0)$, and the map $\sigma_\beta$ has constant sign
    on the cell $(\beta_0,\beta)$.
  \end{itemize}
  We call $\C$ \textbf{piecewise trivial for $\xi$} if $\C^*$ is
  almost piecewise trivial for $\xi^*$.  
\end{df}

\begin{expl}
  \label{radial}
  Let $\zeta := -y \frac{\partial}{\partial x} + x
  \frac{\partial}{\partial y}$, and let $\C$ be the cell decomposition
  of $\RR^2$ consisting of the sets of the form $\{(x,y):\ x \ast 0,\
  y \star 0\}$ with $\ast,\star \in \{=,<,>\}$.  Then $\C$ is
  piecewise trivial for $\zeta$.
\end{expl}

\begin{nrmks}  
  \label{pi-tri_rmks}
  \begin{enumerate}
  \item Any piecewise trivial decomposition for $\xi$ is a Rolle
    decomposition for $\xi$.
  \item If $U$ is bounded, then $\C$ is almost piecewise trivial for
    $\xi$ if and only if $\C$ is piecewise trivial for $\xi$.
  \item We obtain a piecewise trivial decomposition for $\xi$ in the
    following way: first, obtain a $C^p$-cell decomposition $\C$
    compatible with $U$, $\bd(U)$ and $S(\xi)$ satisfying (I) and
    (II).  Then, to satisfy (III)--(V), we only need to refine
    $\Pi_1(\C):= \set{\Pi_1(C):\ C \in \C}$.
  \end{enumerate}
\end{nrmks}

We now fix a piecewise trivial decomposition $\C$ of $\RR^2$ for
$\xi$.  The name ``piecewise trivial'' is justified by:

\begin{thm}  
  \label{piecewise_trivial}
  Let $C \in \C_U$ be open such that $C \cap S(\xi) = \emptyset$.
  Then the relation $\ll_C$ on $\la(C)$ is a linear ordering.
\end{thm}

To prove the theorem, we fix a \textit{bounded}, open $C \in \C_U$
such that $C \cap S(\xi) = \emptyset$.  Establishing the theorem for
this $C$ suffices: if the theorem holds for every bounded, open $D \in
\C$ such that $D \cap S(\xi) = \emptyset$, then the theorem holds with
$\C^*$ and $\xi^*$ in place of $\C$ and $\xi$ (because every $D \in
\C^*$ is bounded).  Since $\phi_2$ is an analytic diffeomorphism, it
follows that the theorem holds for every open $D \in \C$ such that $D
\cap S(\xi) = \emptyset$.

We need quite a bit of preliminary work (see the end of this section
for the proof of the theorem).  For Lemma \ref{H-limit_2} and
Corollary \ref{strip_frontier} below, we fix a $C^p$-curve
$\gamma:[0,1] \into C$ transverse to $\xi$.

\begin{lemma}  
  \label{H-limit_2}
  Let $t_i \in (0,1)$, for $i \in \NN$, such that $t_i \to t \in
  [0,1]$.  Then $C \cap \lim \cl(L(t_i)) = L(t)$.
\end{lemma}

\begin{proof}
  From Proposition \ref{H-limit} we know that $C \cap K$ is a union of
  leaves of $\xi\rest C$, where $K:= \lim\cl(L(t_i))$.  Thus, since
  $\gamma(t_i) \to \gamma(t)$ and $\gamma(t) \in L(t)$, it follows
  that $L(t) \subseteq C \cap K$.  To prove the opposite inclusion, we
  may assume by Proposition \ref{constant_ranks} that every leaf of
  $\xi\rest C$ is the graph of a strictly increasing function (the
  other case is handled similarly).  By Proposition \ref{H-limit}
  again, $\Pi_1(K) = [a,b]$ with $c \le a < b \le d$, and there is a
  continuous function $f:(a,b) \into \RR$ such that $K \cap \big((a,b)
  \times \RR\big) = \gr(f)$.

  Assume for a contradiction that there is a leaf $M$ of $\xi\rest C$
  such that $M \ne L(t)$ and $M \subseteq C \cap K$.  Then $L(t)$ and
  $M$ are disjoint subsets of $\gr(f)$; say $L(t) = \gr(f_t)$, where
  $f_t:(a(t),b(t)) \into \RR$, and $M = \gr(g)$, where $g:(a',b')
  \into \RR$.  We assume here that $a' < b' \le a(t) < b(t)$; the
  other case is again handled similarly.  By our assumption, $c <
  a(t)$ and hence $\lim_{x \to a(t)^+} f_t(x) \in \{\alpha(a(t)),
  \beta(a(t))\}$.  We assume here $\lim_{x \to a(t)^+} f_t(x) =
  \alpha(a(t))$, the other case being handled similarly.  Then by the
  Mean Value Theorem, for every $\epsilon > 0$ there is an $x \in
  (a(t),a(t)+\epsilon)$ such that $f_t'(x) > \alpha'(x)$, that is,
  $\sigma_\alpha(x,f_t(x)) < 0$.  It follows from (V) that
  \begin{itemize}
  \item[($\ast$)] the map $\sigma_\alpha$ has constant negative sign
    on $(\alpha,\alpha_0)$.
  \end{itemize}
  On the other hand, $b' < d$, and we may assume that $\lim_{x \to
    b'^-} g(x) = \alpha(b')$: otherwise, $\lim_{x \to b'^-} g(x) =
  \beta(b')$, and since $$\lim_{x \to a(t)} f(x) = \lim_{x \to a(t)^+}
  f_t(x) = \alpha(a(t)),$$ we can replace $M$ by a leaf of $\xi\rest C$ that is contained in $\gr(f)$ and has the desired property.  But
  $\lim_{x \to b'^-} g(x) = \alpha(b')$ means (as above) that for
  every $\epsilon > 0$ there is an $x \in (b'-\epsilon, b')$ such that
  $g'(x) < \alpha'(x)$, that is, $\sigma_\alpha(x,g(x)) > 0$.  This
  contradicts ($\ast$), so the lemma is proved.
\end{proof}

Put $F:= F(\gamma((0,1)))$; note that $F$ is open by Lemma
\ref{order_along_curve}(1).

\begin{cor}  
  \label{strip_frontier}
  $C \cap \bd(F) = L(0) \cup L(1)$; in particular, there are distinct
  $j_0, j_1 \in \{1,2\}$ such that $C \setminus \cl(F) = U_{L(0),j_0}
  \cup U_{L(1),j_1}$.
\end{cor}

\begin{proof}
  Let $z \in \cl(F) \cap C$, and let $z_i \in F$ be such that $z_i \to
  z$.  Let $t_i \in (0,1)$ be such that $z_i \in L(t_i)$; passing to a
  subsequence if necessary, we may assume that $t_i \to t \in [0,1]$.
  Then $z \in C \cap \lim\cl(L(t_i))$, so $z \in L(t)$ by Lemma
  \ref{H-limit_2}.  Since $F$ is open by Lemma
  \ref{order_along_curve}(1), it follows that $C \cap \bd(F) \subseteq
  L(0) \cup L(1)$.  On the other hand, by Lemma
  \ref{order_along_curve}(2), there is a $j \in \{1,2\}$ such that
  $L(t) \subseteq U_{L(0),j}$ for all $t \in (0,1]$ and $L(t)
  \subseteq U_{1,j'}$ for all $t \in [0,1)$, where $j' \in \{1,2\}
  \setminus \{j\}$.  Hence $L(0) \cup L(1) \subseteq C \cap
  \bd(F(\Gamma))$, and the corollary is proved.
\end{proof}

\begin{df}  
  \label{piecewise_monotone}
  Let $\tau:[0,1] \into U$ be continuous.  We call $\tau$
  \textbf{piecewise $C^p$-monotone in $\xi$} if there are $t_0:= 0 <
  t_1 < t_2 < \cdots < t_k < t_{k+1} := 1$ and $\ast \in \{<,>\}$ such
  that for all $i = 0, \dots, k$, the restriction
  $\tau\rest{(t_i,t_{i+1})}$ is $C^p$, and either $\xi^\perp(\tau(t))
  \cdot \tau'(t) = 0$ for all $t \in (t_i,t_{i+1})$ or
  $\xi^\perp(\tau(t)) \cdot \tau'(t) \ast 0$ for all $t \in
  (t_i,t_{i+1})$.  In this situation, we also say that $\tau$ is
  \textbf{$\ast$-piecewise $C^p$-monotone in $\xi$}.  We call such a
  $\tau$ \textbf{tangent to $\xi$} if each $\tau\rest{(t_i,t_{i+1})}$
  is tangent to $\xi$.
\end{df}

\begin{lemma}  
  \label{piecewise_path_exists}
  Let $v,w \in C$.  Then there is a curve $\tau:[0,1] \into C$ that is
  piecewise $C^p$-monotone in $\xi$ and satisfies $\tau(0) = v$ and
  $\tau(1) = w$.
\end{lemma}

\begin{proof}
  If $L_v = L_w$, then there is a $C^p$-curve $\tau:[0,1] \into L_v$
  such that $\tau(0) = v$ and $\tau(1) = w$, and we are done.  So we
  assume from now on that $L_v \ne L_w$.  Let $j_{vw} \in \{1,2\}$ be
  such that $w \in U_{L_v,j_{vw}}$, and put
  \begin{equation*}
    \ast_{vw} := \begin{cases} < &\text{if } j_{vw} = 1, \\ >
      &\text{if } j_{vw} = 2. \end{cases}
  \end{equation*}
  By o-minimality, there is a definable $C^p$-curve $\tau:[0,1] \into
  C$ such that
  \begin{itemize}
  \item[(I)] $\tau(0) = v$ and $\tau(1) = w$.
  \end{itemize}
  Again by o-minimality, there are $t_0:= 0 < t_1 < \cdots < t_k <
  t_{k+1} := 1$ such that for each $i = 0, \dots, k$,
  \begin{itemize}
  \item[(II)] the map $t \mapsto \xi^\perp(\tau(t)) \cdot \tau'(t)$
    has constant sign on $(t_i,t_{i+1})$.
  \end{itemize}
  By Khovanskii theory \cite{spe:pfaffian}, we may also assume that
  for every $i = 0, \dots, k$,
  \begin{itemize}
  \item[(III)] either $\tau((t_i,t_{i+1})) \cap (L_v \cup L_w) =
    \emptyset$ or $\tau((t_i,t_{i+1})) \subseteq L_v \cup L_w$.
  \end{itemize}
  We now proceed by induction on $k$, simultaneously for all $v,w \in
  C$ and $\tau$ satisfying (I)--(III), to prove that $\tau$ can be
  changed into a curve that is $\ast_{vw}$-piecewise $C^p$-monotone in
  $\xi$.  If $k=0$, then $\tau$ is $\ast_{vw}$-piecewise
  $C^p$-monotone in $\xi$, so we are done.  Therefore, we assume
  that $k>0$ and that the claim holds for lower values of $k$.

  Since $\tau(1) = w \notin L_v$ and $L_v$ is closed in $C$, there is
  a maximal $t \in [0,1)$ such that $\tau(t) \in L_v$, and by our
  choice of $t_1, \dots, t_k$, we have $t = t_i$ for some $i \in \{0,
  \dots, k\}$.  If $i>1$, we replace $\tau\rest{[0,t_i]}$ by a $C^p$
  curve $\tau_1:[0,t_i] \into L_v$ such that $\tau_1(0) = v$ and
  $\tau_1(t_i) = \tau(t_i)$, and we reindex $t_i, \dots, t_{k+1}$ as
  $t_1, \dots, t_{k-i+2}$.  Hence by the inductive hypothesis, we may
  assume that $i \le 1$ and $\tau([0,1]) \subseteq L_v \cup
  U_{L_v,j_{vw}}$.  Put $v':= \tau(t_1)$; we now distinguish two
  cases: \smallskip

  \noindent\textbf{Case 1:} $v' \in L_v$.  Then $\ast_{v'w} =
  \ast_{vw}$, so by the inductive hypothesis (and rescaling), there is
  a curve $\tau_1:[t_1,1] \into C$ that is $\ast_{vw}$-piecewise
  $C^p$-monotone in $\xi$ and satisfies $\tau_1(t_1) = v'$ and
  $\tau_1(1) = w$.  Now replace $\tau\rest{[t_1,1]}$ by $\tau_1$.
  \smallskip

  \noindent\textbf{Case 2:} $v' \notin L_v$.  Then we must have
  $\xi^\perp(\tau(t)) \cdot \tau(t) \ast_{vw} 0$ for all $t \in
  (0,t_1)$.  If $v' \in L_w$, the lemma follows by a similar argument
  as in Case 1, so we assume that $v' \notin L_w$.  We claim again
  that $\ast_{v'w} = \ast_{vw}$ in this situation, from which the
  lemma then follows from the inductive hypothesis as in Case 1.

  To see the claim, by Corollary \ref{strip_frontier} the
  complement of $F(\tau([0,t_1]))$ in $C$ has two connected components
  $U_{L_v,j}$ and $U_{L_{v',j'}}$, where $j,j' \in \{1,2\}$ are
  distinct.  By the above, $j$ must be different from $j_{vw}$, so $w
  \in U_{L_{v',j'}}$, that is, $j' = j_{v'w}$, which implies $j_{vw} =
  j_{v'w}$ as required.
\end{proof}

\begin{lemma}  
  \label{path_replace}
  Let $\tau:[0,1] \into C$ be piecewise $C^p$-monotone in $\xi$
  such that $\tau$ is not tangent to $\xi$.  Then there is a $C^p$
  curve $\gamma:[0,1] \into C$ such that $\gamma$ is transverse to
  $C$, $\gamma(0) = \tau(0)$ and $\gamma(1) = \tau(1)$.
\end{lemma}

\begin{proof}
  Let $t_0:= 0 < t_1 < t_2 < \cdots < t_k < t_{k+1} := 1$ be as in
  Definition \ref{piecewise_monotone}.  We work by induction on $k$;
  if $k=0$, then by hypothesis $\tau$ is transverse to $\xi$, and
  we take $\gamma:= \tau$.  So we assume that $k>0$; for the inductive
  step, it suffices to consider the the case $k=1$.  The hypothesis on
  $\tau$ then implies that at least one of $\tau\rest{(0,t_1)}$ and
  $\tau\rest{(t_1,1)}$ is transverse to $\xi$; so we distinguish
  three cases: \smallskip

  \noindent\textbf{Case 1:} both $\tau\rest{(0,t_1)}$ and
  $\tau\rest{(t_1,1)}$ are transverse to $\xi$.  By Picard's theorem,
  there are an open neighborhood $W \subseteq C$ of $\tau(t_1)$ and a
  $C^p$-diffeomorphism $f:\RR^2 \into W$ such that $f(0) = \tau(t_1)$
  and $f^*\xi = \partial/\partial x$, where $f^*\xi$ is the pull-back
  of $\xi$ via $f$.  Then for some $\epsilon > 0$, the continuous
  curve $f^{-1} \circ \tau \rest{(t_1-\epsilon,t_1+\epsilon)}$ is
  $C^p$ and transverse to $\partial/\partial x$ on $(t_1-\epsilon,t_1)
  \cup (t_1, t_1+\epsilon)$.  Using standard smoothing arguments from
  analysis, we can now find a $C^p$-curve
  $\eta:(t_1-\epsilon,t_1+\epsilon) \into \RR^2$ that is transverse to
  $\partial/\partial x$ and satisfies $\eta(t) = f^{-1}(\tau(t))$ for
  all $t \in (t_1-\epsilon,t_1-\epsilon/2) \cup
  (t_1+\epsilon/2,t_1+\epsilon)$.  Now define $\gamma:[0,1] \into C$
  by
  \begin{equation*} \gamma(t):= 
    \begin{cases} 
      \tau(t) &\text{if } 0
      \le t < t_1-\epsilon \text{ or } t_1 + \epsilon < t \le 1, \\
      f(\eta(t)) &\text{if } t_1-\epsilon \le t \le
      t_1+\epsilon. 
    \end{cases}
  \end{equation*}
  \smallskip

  \noindent\textbf{Case 2:} $\tau\rest{(0,t_1)}$ is transverse to
  $\xi$ and $\tau\rest{(t_1,1)}$ is tangent to $\xi$.  Since
  $\tau([t_1,1])$ is compact, there are (by Picard's theorem again)
  $s_0:= t_1 < s_1 < \cdots < s_l < s_{l+1}:= 1$, open neighborhoods
  $W_i \subseteq U$ of $\tau(s_i)$ and $C^p$-diffeomorphisms
  $f_i:\RR^2 \into W_i$, for $i=0, \dots, l+1$, such that
  $\tau([t_1,1]) \subseteq W_0 \cup \cdots \cup W_{l+1}$, $f_i(0) =
  \tau(s_i)$ and $f_i^*\xi = \partial/\partial x$ for each $i$.  We
  assume that $l=0$, so that $s_0 = t_1$ and $s_1 = 1$; the general
  case then follows by induction on $l$.

  Let $u \in (t_1,1)$ be such that $\tau(u) \in W_0 \cap W_1$.
  Working with $f_0$ similarly as in Case 1, we can replace
  $\tau\rest{[0,u]}$ by a $C^p$-curve $\eta:[0,u] \into C$ transverse
  to $\xi$ such that $\eta(0) = \tau(0)$ and $\eta(u) = \tau(u)$.
  Define $\eta(t):= \tau(t)$ for $t \in (u,1]$; repeating the
  procedure with $\eta$ and $f_1$ in place of $\tau$ and $f_0$, we
  obtain a $C^p$-curve $\gamma:[0,1] \into C$ that is transverse to
  $\xi$ and satisfies $\gamma(0) = \tau(0)$ and $\gamma(1) =
  \tau(1)$, as desired.  \smallskip

  \noindent\textbf{Case 3:} $\tau\rest{(0,t_1)}$ is tangent to
  $\xi$ and $\tau\rest{(t_1,1)}$ is transverse to $\xi$.  This
  case is similar to Case 2.
\end{proof}

%\begin{lemma}  \label{equal_foliation}
% Let $\tau:[0,1] \into C$ be piecewise $C^p$-monotone in $\xi$, and
% let $\gamma:[0,1] \into C$ be $C^p$ and transverse to $\xi$ such
% that $\gamma(0) = \tau(0)$ and $\gamma(1) = \tau(1)$.  Then
% $F(\Gamma) = F(\T)$, where $\T:= \tau([0,1])$.
%\end{lemma}
%
%\begin{proof}
% Note that the sets $\tau^{-1}(L)$ and $\gamma^{-1}(L)$ are connected
% for every $L \in \la(C)$.  We show that $F(\Gamma) \subseteq F(\T)$;
% the other inclusion is proved the same way.
%% Let $z \in F(\Gamma)$; we need to show that $L_z \cap \T \ne
% \emptyset$.  This is clear if $z \in L(0)$ or $z \in L(1)$, so we
% assume that $z \in L(t)$ for some $t \in (0,1)$.  Since
% $\gamma^{-1}(L(t))$ is connected, the points $\tau(0) = \gamma(0)$
% and $\tau(1) = \gamma(1)$ belong to distinct components of $C
% \setminus L(t)$.  Hence there is a $t \in (0,1)$ such that $\tau(t)
% \in L(t) = L_z$.
%\end{proof}

Combining Lemmas \ref{piecewise_path_exists} and \ref{path_replace},
we obtain:

\begin{cor}  
  \label{transverse_path}
  Let $u,v \in C$ be such that $L_u \ne L_v$.  Then there is a $C^p$
  curve $\gamma:[0,1] \into C$ such that $\gamma(0) = u$, $\gamma(1) =
  v$ and $\gamma$ is transverse to $\xi$. \qed
\end{cor}

\begin{proof}[Proof of Theorem \ref{piecewise_trivial}]
  Let $M,L \in \la(C)$ be distinct and choose $v \in M$ and $w \in L$.
  By Corollary \ref{transverse_path}, there is a $C^p$-curve
  $\gamma:[0,1] \into C$ such that $\gamma(0) = v$, $\gamma(1) = w$
  and $\gamma$ is transverse to $\xi$.  Hence $t \mapsto
  \xi^\perp(\gamma(t)) \cdot \gamma'(t)$ has constant nonzero sign on
  $[0,1]$; this shows that $\ll_C$ is irreflexive.  Transitivity
  follows by a similar argument.
\end{proof}

\section{Foliation orderings}  \label{foliation_orderings}

Let $\xi = a_1\frac{\partial}{\partial x} + a_2
\frac{\partial}{\partial y}$ be a definable vector field of class
$C^1$ on $\RR^2$.  We fix a piecewise trivial decomposition $\C$ of
$\RR^2$ for $\xi$; refining $\C$ if necessary, we may assume that $\C$
is a stratification.  To simplify statements, we put
\begin{equation*}
  \Creg := \set{C \in \C:\ C \cap S(\xi) = \emptyset}.
\end{equation*}
For instance in Example \ref{radial}, the piecewise trivial
decomposition $\C$ is a stratification and $\Creg = \C \setminus
\{0\}$.

\begin{nrmk}  
  \label{strat_rmk}
  $\C$ being a stratification has the following consequence: for every
  $1$-dimensional $C \in \C$, there are exactly two distinct open $D
  \in \C$ such that $C \cap \fr(D) \ne \emptyset$, and for each of
  these $D$ we have $C \subseteq \fr(D)$.
\end{nrmk}

Let $V \subseteq \RR^2 \setminus S(\xi)$ be an integral manifold of
$\xi$, that is, a $1$-dimensional manifold tangent to $\xi$.  Given
$u,v \in V$, we define $u <_V^\xi v$ if and only if there is a $C^1$
path $\gamma:[0,1] \into V$ such that $\gamma(0) = u$, $\gamma(1) = v$
and $\xi(\gamma(t)) \cdot \gamma'(t) > 0$ for all $t \in [0,1]$.

\begin{lemma}  
  \label{dlo_1}
  Assume that $V$ is connected and not a compact leaf.  Then the
  relation $<_V^\xi$ defines a dense linear ordering of $\,V$
  without endpoints.
\end{lemma}

\begin{proof}
  Let $u,v \in V$ be such that $u \ne v$.  Since $V$ is connected, we
  get $u <^\xi_V v$ or $v <^\xi_V u$.  On the other hand, if there are
  $C^1$-paths $\gamma, \delta:[0,1] \into V$ such that $\gamma(0) =
  \delta(1) = u$, $\gamma(1) = \delta(0) = v$ and $\xi(\gamma(t))
  \cdot \gamma'(t) > 0$ and $\xi(\delta(t)) \cdot \delta'(t) > 0$ for
  all $t \in [0,1]$, then $\gamma([0,1]) \cup \delta([0,1])$ is a
  compact leaf of $\xi$ contained in $V$; since $V$ is connected, it
  follows that $V$ is a compact leaf, a contradiction.
\end{proof}

We now fix a $C \in \Creg$ such that $\dim(C) > 0$.

\begin{df}  
  \label{ordering}
  The foliation of $\xi$ induces an ordering $<_C^\xi$ on $C$ as
  follows:
  \begin{itemize}
  \item Suppose that $C$ is open, and let $u,v \in C$.  Then every
    leaf of $\xi\rest{C}$ is non-compact by Proposition
    \ref{constant_ranks}.  Thus, we define $u <_C^\xi v$ if and only
    if $L_u \ll_C^\xi L_v$ or $L_u = L_v$ and $u <_{L_u}^\xi v$.
  \item Suppose that $\dim(C) = 1$ and $C$ is tangent to $\xi$.
    Then $C$ is a connected, non-compact integral manifold of
    $\xi$, so we define $<_C^\xi$ as before Lemma \ref{dlo_1}.
  \item Suppose that $\dim(C) = 1$ and $C$ is transverse to $\xi$.
    Let $u,v \in C$; we define $u <^\xi_C v$ if and only if there is a
    $C^1$-curve $\gamma:[0,1] \into C$ such that $\xi^\perp(\gamma(t))
    \cdot \gamma'(t) > 0$ for all $t \in [0,1]$.
  \end{itemize}
  As before, we omit the superscript $\xi$ whenever it is clear
  from context.

  A \textbf{$<_C$-interval} is a set $A$ of the form $(a,b) := \set{c
    \in C:\ a \ast_1 c \ast_2 b}$ with $a,b \in C$, or $(a,\infty) :=
  \set{c \in C:\ a \ast c}$ with $a \in C$, or $(-\infty,b) := \set{c
    \in C:\ c \ast c}$ with $b \in C$, where $\ast, \ast_1, \ast_2 \in
  \{<_C, \le_C\}$; we call $A$ \textbf{open} if $\ast = \ast_1 =
  \ast_2 = <_C$.
\end{df}

\begin{lemma}  
  \label{dlo_2}
  The ordering $<_C$ is a dense linear ordering on $C$ without
  endpoints.  Moreover, if $\dim(C) = 1$, then every $<_C$-bounded
  subset of $C$ has a least upper bound.
\end{lemma}

\begin{proof}
  It is clear from the definition that $C$ has no endpoints with
  respect to $<_C$.  Density and linearity follow from Lemmas
  \ref{order_along_curve} and \ref{dlo_1} if $\dim(C) = 1$, and if $C$
  is open, they follow from Lemma \ref{dlo_1} and Theorem
  \ref{piecewise_trivial}.

  For the second statement, assume that $\dim(C) = 1$ and let
  $\alpha:(0,1) \into \RR^2$ be $C^1$ and injective such that $C =
  \alpha((0,1))$.  If $C$ is tangent to $\xi$, then the map $t \mapsto
  \xi(\alpha(t)) \cdot \alpha'(t)$ has constant nonzero sign, and if
  $C$ is transverse to $\xi$, then the map $t \mapsto
  \xi^\perp(\alpha(t)) \cdot \alpha'(t)$ has constant nonzero sign.
  Thus in both cases, the map $\alpha:\big((0,1),<\big) \into (C,<_C)$
  is either order-preserving or order-reversing; the second statement
  follows.
\end{proof}

We assume for the remainder of this section that either $C$ is open,
or $C$ is $1$-dimensional and tangent to $\xi$.

\begin{df}  
  \label{projections}
  For each leaf $L$ of $\xi\rest{C}$, it follows from Proposition
  \ref{constant_ranks} that $\fr(L)$ consists of exactly two points
  $P_L^>, P_L^< \in \fr(C) \cup \{\infty\}$, where, for $\ast \in
  \{>,<\}$, $P_L^\ast$ is the unique of these two points with the
  property that for every $C^1$-curve $\gamma:[0,1) \into L$
  satisfying $\gamma(0) \in L$ and $\lim_{t \to 1} \gamma(t) =
  P_L^\ast$, we have $\xi(\gamma(t)) \cdot \gamma'(t) \ast 0$
  for all $t \in [0,1)$.  In this situation, we define the
  \textbf{forward projection} $\f_C:C \into \fr(C) \cup \{\infty\}$
  and the \textbf{backward projection} $\b_C:C \into \fr(C) \cup
  \{\infty\}$ as
  \begin{equation*}
    \f_C(z) := P_{L_z}^> \quad\text{and}\quad \b_C(z) := P_{L_z}^<,
    \quad\text{for all } z \in C.
  \end{equation*}
\end{df}

From now on we assume that $C$ is open, and we let $D \in \Creg$ be of
dimension $1$ and contained in $\fr(C)$ such that $D$ is transverse to
$\xi$.

\begin{lemma}  
  \label{projection_image}
  Either $D \subseteq \f_C(C)$ and $D \cap \b_C(C) = \emptyset$, or $D
  \subseteq \b_C(C)$ and $D \cap \f_C(C) = \emptyset$.
\end{lemma}

\begin{proof}
  Let $\alpha:(0,1) \into \RR^2$ be a definable $C^1$-map such that $D
  = \alpha((0,1))$ and $\xi^\perp(\alpha(t)) \cdot \alpha'(t) > 0$ for
  all $t \in (0,1)$.  Thus, either $\xi(\alpha(t))$ points into $C$
  for all $t$, or $\xi(\alpha(t))$ points out of $C$ for all $t$.  In
  the first case, we have $\f_C(C) \cap D = \emptyset$, and in the
  second case $\b_C(C) \cap D = \emptyset$.  Moreover by Picard' s
  Theorem, for every $w \in D$ there is an integral manifold $V
  \subseteq \RR^2$ of $\xi$ such that $V \cap D = \{w\}$; hence,
  either $w \in \f_C(C)$ or $w \in \b_C(C)$.
\end{proof}

\begin{lemma}  
  \label{increasing_projections}
  The maps $\f_C\rest{\f_C^{-1}(D)}$ and $\b_C\rest{\b_C^{-1}(D)}$ are
  increasing.
\end{lemma}

\begin{proof}
  We prove the lemma for $\f_C$.  Let $u,v \in C$ with $u <_C v$ be
  such that $\f_C(u), \f_C(v) \in D$; we may clearly assume that $L_u
  \ll_C L_v$, and hence (by Picard's Theorem) that $\f_C(u) \ne
  \f_C(v)$.

  We assume here that $D = \gr(\alpha)$, where $\alpha:(a,b) \into
  \RR$ is a definable $C^1$-function; the case $D = \{a\} \times
  (b,c)$ is handled similarly.  Let also $\beta:(a,b) \into \RR$ be a
  definable $C^1$-function such that $C = (\alpha,\beta)$ or $C =
  (\beta,\alpha)$; we assume here the former, the latter being handled
  similarly.  For $s \in [0,1]$, we put
  \begin{equation*}
    \alpha_s(t):= (1-s) \alpha(t) + s \beta(t), \quad a < t < b.
  \end{equation*}
  Then for every $t \in (a,b)$, we have $\lim_{s \to 0} \alpha_s(t) =
  \alpha(t)$ and $\lim_{s \to 0} \alpha_s'(t) = \alpha'(t)$.

  Let now $a < a' < b' < b$ be such that $\f_C(u), \f_C(v) \in
  \gr\alpha\rest{(a',b')}$.  Since $D$ is transverse to $\xi$,
  there is an $\epsilon > 0$ such that $\gr\alpha_s\rest{(a',b')}$ is
  transverse to $\xi$ for all $s \in [0,\epsilon)$.  It follows
  from the previous paragraph that the map $t \mapsto
  \sigma_\alpha(t,\alpha(t))$ has the same constant nonzero
  sign as the map $t \mapsto \sigma_{\alpha_s}(t,\alpha_s(t))$,
  for all $s \in (0,\epsilon)$.  Therefore by Lemma
  \ref{order_along_curve}(2) and the definition of $<_D$, we have
  $\f_C(u) <_D \f_C(v)$, as required.
\end{proof} 

\begin{cor}  
  \label{projection_of_interval}
  Let $I \subseteq C$ be a $<_C$-interval.  Then each of $\f_C(I) \cap
  D$ and $\b_C(I) \cap D$ is either empty, a point or an open
  $<_D$-interval.
\end{cor}

\begin{proof}
  Assume that $a,b \in \f_C(I) \cap D$ are such that $a <_D b$, and
  let $c \in D$ be such that $a <_D c <_D b$; it suffices to show that
  $c \in \f_C(I)$.  By Lemma \ref{projection_image}, $c \in \f_C(C)$.
  Let $u, v, w \in C$ be such that $a = \f_C(u)$, $b = \f_C(v)$, $c =
  \f_C(w)$ and $u,v \in I$.  Then $u <_C w <_C v$ by Lemma
  \ref{increasing_projections}, as required.
\end{proof} 

We fix a set $E_C \subseteq C$ such that $|E_C \cap L| = 1$
for every $L \in \la(C)$ and put $<_{E_C}:= <_C\rest{E_C}$, and we
denote by $e_L$ the unique element of $E \cap L$, for every $L \in
\la(C)$.

\begin{rmk}
  The map $L \mapsto L \cap E_C : (\la(C),\ll_C) \into (E_C,<_{E_C})$
  is an isomorphism of ordered structures.
\end{rmk}

\begin{prop}  
  \label{inverse_image}
  Let $\g \in \{\f,\b\}$.  If $D \subseteq \g_C(C)$, then $D_\g:=
  \g_C^{-1}(D) \cap E_C$ is an $<_{E_C}$-interval, and the map
  $\g_C\rest{D_\g}: (D_\g,<_{E_C}\rest{D_\g}) \into (D,<_D)$ is an
  isomorphism of ordered structures.
\end{prop}

\begin{proof}
  The transversality of $D$ to $\xi$ implies that if $u \in D$ and
  $L_1,L_2 \in \la(C)$ are such that $u = P^>_{L_1} = P^>_{L_2}$ or $u
  = P^<_{L_1} = P^<_{L_2}$, then $L_1 = L_2$.  Thus by Lemma
  \ref{increasing_projections}, the map $\g_C\rest{D_\f}$ is strictly
  increasing, so the lemma follows.
\end{proof}

\section{Progression map}  \label{progression_map}

We continue working with $\xi$ and $\C$ as in Section
\ref{foliation_orderings}, and we adopt all corresponding notations.
We let
\begin{renumerate}
\item $\Copen$ be the collection of all open cells in $\Creg$;
\item $\Ctan$ be the collection of all cells in $\Creg$ that are of
  dimension $1$ and tangent to $\xi$;
\item $\Ctrans$ be the collection of all cells in $\Creg$ that are of
  dimension $1$ and transverse to $\xi$; and
\item $\Csingle$ the collection of all $p \in \RR^2$ such that $\{p\}
  \in \Creg$.
\end{renumerate}

By Lemma \ref{projection_image} and since $\C$ is a stratification,
there are, for each $C \in \Ctrans$, distinct and unique cells $C^\b,
C^\f \in \Copen$ such that $C \cap \cl(C^\b) \ne \emptyset$, $C \cap
\cl(C^\f) \ne \emptyset$ and
\begin{equation*}
  C \subseteq \f_{C^\b}(C^\b) \text{ and } C \subseteq
  \b_{C^\f}(C^\f).
\end{equation*}
Similarly, there are, for each $p \in \Csingle$, distinct and unique
cells $p^\b, p^\f \in \Copen \cup \Ctan$ such that $p \in \cl(p^\b)$, $p
\in \cl(p^\f)$ and
\begin{equation*}
  p \in \f_{p^\b}(p^\b) \text{ and } p \in \b_{p^\f}(p^\f).
\end{equation*}
(For $p \in \Csingle$, we use the fact that there is an open box $B$
containing $p$ such that the leaf of $\xi\rest{B}$ passing through
$p$ is a Rolle leaf.)  For each $C \in \Ctan$, we fix an arbitrary
element $e_C \in C$; note that for each $z \in C$, $C$ is the unique
leaf $L_z$ of $\xi\rest{C}$ containing $z$.

We now define $\f',\b':\RR^2 \into \RR^2 \cup \{\infty\}$ by
\begin{equation*}
  \f'(z):= \begin{cases} \f_{C}(z) &\text{if } z \in C \in \Copen \cup
    \Ctan \text{ and } e_{L_z} \le_{L_z} z, 
    \\ e_{L_z} &\text{if } z \in C \in \Copen \cup
    \Ctan \text{ and } z <_{L_z} e_{L_z}, \\
    \big(\b_{C^\f}\rest{E_{C^\f}}\big)^{-1}(z) &\text{if } z \in C
    \in \Ctrans \cup \Csingle, \\ z &\text{if }
    z \in S(\xi) \end{cases}
\end{equation*}
and
\begin{equation*}
  \b'(z):= \begin{cases} \b_{C}(z) &\text{if } z \in C \in \Copen \cup
    \Ctan \text{ and } z \le_{L_z} e_{L_z}, 
    \\ e_{L_z} &\text{if } z \in C \in \Copen \cup
    \Ctan \text{ and } e_{L_z} <_{L_z} z, \\
    \big(\f_{C^\b}\rest{E_{C^\b}}\big)^{-1}(z) &\text{if } z \in C
    \in \Ctrans \cup \Csingle, \\ z &\text{if }
    z \in S(\xi). \end{cases}
\end{equation*}

\begin{df}
  \label{progression_def}
  We define $\f,\b:\RR^2\cup\{\infty\} \into \RR^2\cup\{\infty\}$ by
  \begin{equation*}
    \f(z):= \begin{cases} \f'(z) &\text{if } z \in \RR^2 \text{ and }
      \f'(z) \notin S(\xi), \\  \infty &\text{otherwise} \end{cases}
  \end{equation*}
  and
  \begin{equation*}
    \b(z):= \begin{cases} \b'(z) &\text{if } z \in \RR^2 \text{ and }
      \b'(z) \notin S(\xi), \\ \infty &\text{otherwise}. \end{cases}
  \end{equation*}
  We call $\f$ a \textbf{progression map} associated to $\xi$ and
  $\b$ a \textbf{reverse progression map} associated to $\xi$.  We
  put
  \begin{equation*}
    \mathcal{C}_1= \Ctrans \cup \Csingle \cup \bigcup\set{E_C:\ C \in
      \Copen} \cup \set{\{e_C\}:\ C \in \Ctan}
  \end{equation*}
  and let $B:=\bigcup \mathcal{C}_1$; note that $\f(\RR^2) \subseteq B
  \cup \{\infty\}$ and $\b(\RR^2) \subseteq B \cup \{\infty\}$.
  Finally, we define $\f^0:\RR^2\cup\{\infty\} \into
  \RR^2\cup\{\infty\}$ by $\f^0(x):= x$, and for $k > 0$ we define
  $\f^k:\RR^2\cup\{\infty\} \into \RR^2\cup\{\infty\}$ inductively on
  $k$ by $\f^k(x) := \f(\f^{k-1}(x))$.
\end{df}

\begin{prop} 
  \label{intersect} 
  Let $X \in \mathcal{C}_1$ and $L$ be a compact leaf of $\xi$.
  Then $|X \cap L| \leq 1$.
\end{prop}

\begin{proof}
  If $X \in \Csingle$ or $X = \{e_C\}$ for some $C \in \Ctan$, the
  conclusion is trivial.  By Lemma \ref{separating}(2), $L$ is a Rolle
  leaf of $\xi$; in particular, $|X \cap L| \le 1$ if $X \in
  \Ctrans$.  So we may assume that $X = E_C$ for some $C \in \Copen$.
  Then there is at most one $L' \in \la(C)$ contained in $L$:
  otherwise by Corollary \ref{transverse_path}, there is a $C^1$-curve
  $\gamma:[0,1] \into C$ transverse to $\xi$ such that $\gamma(0),
  \gamma(1) \in L$, a contradiction.  It follows again that $|X \cap L|
  \le 1$. 
\end{proof}

\begin{prop}  
  \label{fixed_points}
  There is an $N \in \NN$ such that for every $x \in B$, the leaf of
  $\xi$ through $x$ is compact if and only if $\,\f^N(x) = x$.
\end{prop}

\begin{proof}
  Let $x \in B$; if $\f^k(x) = x$ for some $k>0$, then the leaf of
  $\xi$ through $x$ is compact.  For the converse, we assume that
  the leaf $L$ of $\xi$ through $x$ is compact.  Since $L$ is
  compact, we have $L \cap S(\xi) = \emptyset$, that is, $\f^k(x)
  \in B$ for every $k > 0$.  Thus with $n := |\Creg|+1$, there are a
  $C \in \Creg$ and $0 \le k_1 < k_2 \le n$ such that $\f^{k_1}(x),
  \f^{k_2}(x) \in C$.  It follows from Proposition \ref{intersect}
  that $\f^{k_1}(x) = \f^{k_2}(x)$, and hence that 
  \begin{equation*}
    x = \b^{k_1} \circ \f^{k_1}(x) = \b^{k_1} \circ \f^{k_2}(x) =
    \f^{k_2-k_1}(x). 
  \end{equation*}
  Since $n$ is independent of $x \in B$, the number $N := n!$ will do.
\end{proof}

\section{Flow configuration theories}  \label{flowconfig}

Inspired by the previous sections, we now define a first-order
theory as described in the introduction.  Our main goal, reached in
Section \ref{QE}, is to show that this theory admits quantifier
elimination in a language suitable to our purposes.

\begin{df}  
  \label{config}
  A \textbf{flow configuration} is a tuple $$\Phi = (\phiopen,
  \phitan, \phitrans, \phisingle, \phi^\b, \phi^\f, \min, \max,
  N_\Phi)$$ such that $\phiopen$, $\phitan$, $\phitrans$ and
  $\phisingle$ are pairwise disjoint, finite sets, 
  \begin{gather*}
    \phi^\b, \phi^\f : \phitrans \cup \phisingle \into \phiopen \cup
    \phitan, \\ \min, \max: \phiopen \cup \phitan \cup \phitrans \into
    \phisingle \cup \{\infty\}
  \end{gather*}
  and $N_\Phi \in \NN$.  In this situation, we shall write $a^\b$ and
  $a^\f$ instead of $\phi^\b(a)$ and $\phi^\f(a)$, for $a \in \phitrans
  \cup \phisingle$.
\end{df}

\begin{expl}
  \label{omega_config}
  Let $\xi$ be a vector field on $\RR^2$ of class $C^1$ and definable
  in an \hbox{o-minimal} expansion of the real field, and let $\C$ be
  a piecewise trivial cell decomposition of $\RR^2$ that is also a
  stratification.  We define $\Copen$, $\Ctan$, $\Ctrans$, $\Csingle$
  and $^\b,\, ^\f : \Ctrans \cup \Csingle \into \Copen \cup \Ctan$ as in
  Section \ref{progression_map}, and we let $N \in \NN$ be as in
  Proposition \ref{fixed_points}.  

  Let $C \in \Copen \cup \Ctan \cup \Ctrans$.  If there is a point in
  $\Csingle$ that is contained in the closure of every set $\set{x \in
    C:\ x <_C^\xi a}$ with $a \in C$, we let $\min(C)$ be any such
  point; otherwise, we put $\min(C) := \infty$.  Similarly, if there
  is a point in $\Csingle$ that is contained in the closure of every
  set $\set{x \in C:\ a <_C^\xi x}$ with $a \in C$, we let
  $\max(C)$ be any such point; otherwise, we put $\max(C) := \infty$.
  Then the tuple $$\Phi_\xi = \Phi_\xi(\C) := (\Copen, \Ctan, \Ctrans,
  \Csingle,\, ^\b,\, ^\f, \min, \max, N)$$ is a \textbf{flow
    configuration associated to $\xi$}.
\end{expl}

For the remainder of this section, we fix a flow configuration $\Phi$.

\begin{df}
  \label{flow_language}
  Let $\la(\Phi)$ be the first-order language consisting of
  \begin{renumerate}
  \item a unary predicate $C$ and a binary predicate $<_C$, for each
    $C \in \phiopen \cup \phitan \cup \phitrans$;
  \item a unary predicate $E_C$ for each $C \in \phiopen$ and a
    constant symbol $e_C$ for each $C \in \phitan$;
  \item a constant symbol $s$, and a constant symbol $c$ for each $c
    \in \phisingle$;
  \item unary function symbols $\f$ and $\b$;
  \item constant symbols $r^\g_C$ and $s^\g_C$ for each $C \in
    \phitrans$ and $\g \in \{\f,\b\}$.
  \end{renumerate}
  Throughout the rest of this paper, for $m \in \NN$ we write $\f^m$
  for the $\la(\Phi)$-word consisting of $m$ repetitions of the symbol
  $\f$, and similarly for $\b^m$.
\end{df}

\begin{expl}
  \label{omega_structure}
  Let $\xi$ and $\C$ be as in Example \ref{omega_config}; we adopt the
  notations used there.  We associate to $\xi$ a unique
  $\la(\Phi_{\xi})$-structure $\M_{\xi} = \M_\xi(\C)$ as follows:
  \begin{renumerate}
  \item the universe $M_\xi$ of $\M_{\xi}$ is $\RR^2 \setminus
    S(\xi) \cup \{\infty\}$;
  \item for each $C \in \Copen \cup \Ctan \cup \Ctrans$, the predicate
    $C$ is interpreted by the corresponding cell in $\C$, and the
    predicate $<_C$ is interpreted by the union of $<_C^\xi$ with
    $\set{(\min(C),a):\ a \in C}$ and  $\set{(a,\max(C)):\ a \in C}$;
  \item for each $C \in \Copen$, the predicate $E_C$ is interpreted by
    the set $E_C$ described in Section \ref{progression_map}, and for
    each $C \in \Ctan$, the constant $e_C$ is interpreted by the
    element $e_C \in C$ picked in Section \ref{progression_map};
  \item the constant $s$ is interpreted as $\infty$, and for each $c
    \in \Csingle$, the constant $c$ is interpreted as the
    corresponding element of $\Csingle$;
  \item the functions $\f$ and $\b$ are interpreted by the
    corresponding forward progression and reverse progression maps;
  \item for each $C \in \Ctrans$ and $\g \in \{\f,\b\}$, the constants
    $r^\g_C$ and $s^\g_C$ are interpreted as the lower and upper
    endpoints, respectively, of the interval $\g(C)$ in $E_{C^\g} \cup
    \{\min(C^\g), \max(C^\g)\}$.
  \end{renumerate}
\end{expl}

\begin{df}
  \label{flow_theory}
  We put $\Phi_0:= \phiopen \cup \phitan \cup \phitrans$; intending to
  capture the theory of the previous example, we let $T(\Phi)$ be the
  $\la(\Phi)$-theory consisting of the universal closures of the
  formulas in the axiom schemes (F1)--(F15) below.
  \begin{lenumerate}{F}{0}
  \item The formulas
    \begin{enumerate}
    \item ${\displaystyle \bigwedge_{c,d \in \phisingle, c \ne d} \neg c=d \wedge
      \bigwedge_{c \in \phisingle, C \in \Phi_0} \neg C(c)}$,
  \item ${\displaystyle \bigwedge_{c \in \phisingle} \neg c=s \wedge
      \bigwedge_{C \in \Phi_0} \neg C(s)}$,
  \item ${\displaystyle x=s \vee \bigvee_{c \in \phisingle} x = c \vee
      \bigvee_{C \in \Phi_0} \left(C(x) \wedge \bigwedge_{D \in \Phi_0,
          D \ne C} \neg D(x) \right)}$.
    \end{enumerate}
  \item For each $C \in \Phi_0$ the sentences stating that $<_C$ is a
    dense linear ordering of $C$, together with $C(x) \rightarrow (x <_C
    \max(C) \wedge \min(C) <_C x)$.
  \end{lenumerate}

  \begin{rmk}
    We do not wish to state that $<_C$ is a linear order on all of $C
    \cup \{\min(C),\max(C)\}$, because it is possible that
    $\min(C)=\max(C)$.  The axioms (F2) suffice for our purpose, which
    is to be able to refer to $C$ as the $<_C$-interval between
    $\min(C)$ and $\max(C)$.
  \end{rmk}

  \begin{lenumerate}{F}{2}
  \item The formula ${\displaystyle \bigwedge_{C \in \phitan} C(e_C) \wedge
    \bigwedge_{C \in \phiopen} E_C(x) \to C(x)}$.
  \item For each $C \in \phiopen$ the sentences stating that the
    restriction of $<_C$ to $E_C$ is a dense linear ordering.
  \item For each $(\g,\h) \in \{(\f,\b), (\b,\f)\}$ and $\ast \in
    \{\le, \ge\}$ the formulas
    \begin{enumerate}
    \item $\g(s) = s \wedge (\neg x=s \to \neg \g(x) = x)$, 
    \item ${\displaystyle \bigwedge_{c \in \phisingle} (\neg \g(c) = s
        \to \h(\g(c)) = c)}$,
    \item ${\displaystyle \bigwedge_{C \in \phiopen} C \g(x) \to E_C(
        \g(x)) \wedge \bigwedge_{C \in \phitan} C(\g(x)) \to
        \g(x) = e_C}$,
    \item ${\displaystyle \bigwedge_{C \in \phitan} (C(x) \wedge e_C
        \ast_C x \ast_C \g(e_C)) \to \g(x) = \g(e_C)}$,
    \item ${\displaystyle \bigwedge_{C \in \phitan} (C(x) \wedge e_C
        \ast_C x \ast_C \h(e_C)) \to \g(x) = e_C}$.
    \end{enumerate}
  \item For each $C \in \Copen$ and $\g \in \{\f,\b\}$ the formula
    \begin{equation*}
      (E_C(x) \wedge E_C(y) \wedge \g(x) = \g (y)) \to (\g (x) = s \vee
      x=y).
    \end{equation*}
  \item For each $c \in \phisingle$ and $\g \in \{\f,\b\}$, the
    sentences $\g (c) = e_{c^\g}$ if $c^\g \in \phitan$ and $E_{c^\g}(
    \g (c))$ if $c^\g \in \phiopen$.
  \item For each $C \in \phitrans$ and $(\g,\h) \in \{(\f,\b),
    (\b,\f)\}$ the sentences stating that $\g(C)$ is an interval $I_1$
    in $E_{C^\g}$ and $\g\rest{C}:C \into I_1$ is an order-isomorphism.
  \item For each $C \in \phiopen$ and $(\g,\h) \in \{(\f,\b),
    (\b,\f)\}$ the formula
    \begin{equation*}
      E_C(x) \to \left(\g (x) = s \vee \bigvee_{D \in
          \phitrans,\ C = D^\h} D (\g (x)) \vee \bigvee_{d \in
          \phisingle,\ C = d^\h} \g (x) = d \right).
    \end{equation*}
  \end{lenumerate}

  We need more axioms describing the ordering $<_C$ and the behavior
  of $\f$ and $\b$ on $C$, for $C \in \phiopen$.  For example, if $x
  \in C \setminus E_C$, we want that $x$ has either a unique
  predecessor or a unique successor in $E_C$.  Also, for any $y \in
  E_C$, the set of points $x$ for which $y$ is either the predecessor
  or successor is infinite and densely ordered by $<_C$.  For
  convenience, we let $\phi_C^\f(x,y)$ be the formula
  \begin{equation*}
    C(x) \wedge \neg E_C(x) \wedge E_C(y)
    \wedge x <_C y \wedge \neg \exists z(E_C(z) \wedge x<_Cz<_Cy)
  \end{equation*}
  and $\phi_C^\b(x,y)$ be the formula
  \begin{equation*}
    C(x) \wedge \neg E_C(x) \wedge E_C(y)
    \wedge y <_C x \wedge \neg \exists z(E_C(z) \wedge y<_Cz<_Cx).
  \end{equation*}

  \begin{lenumerate}{F}{9}
  \item For each $C \in \phiopen$ the formulas
    \begin{enumerate}
    \item ${\displaystyle C(x) \wedge \neg E_C(x) \rightarrow \exists y
        (\phi^\f_C(x,y) \vee \phi_C^\b(x,y))}$,
    \item ${\displaystyle \exists y \phi_C^\f(x,y) \rightarrow \neg
        \exists y \phi_C^\b(x,z)}$,
    \item ${\displaystyle \exists y \phi_C^\b(x,y) \rightarrow \neg
        \exists y \phi_C^\f(x,y)}$,
    \end{enumerate}
    and the formula scheme $E_C(y) \rightarrow
    \exists^{\infty}x\phi_C^\f(x,y) \wedge
    \exists^{\infty}x\phi_C^\b(x,y)$.

  \item For each $C \in \phiopen$ the sentences stating that for every
    $y \in E_C$, the restriction of $<_C$ to the set $C_y := \{x:
    \phi_C^\b(x,y) \vee \phi_C^\f(x,y) \vee x=y\}$ is a dense linear
    ordering, together with $C_y(x) \rightarrow (x <_C \f (y) \wedge
    \g (y) <_C x)$.
  \item For each $C \in \phiopen$ and $(\g,\h) \in \{(\f,\b),
    (\b,\f)\}$ the formulas
    \begin{enumerate}
    \item ${\displaystyle C(x) \wedge \neg E_C(x) \wedge \exists y
        \phi_C^\g(x,y) \rightarrow \forall z(\phi_C^\g(x,z) \rightarrow
        \g (x) =z)}$, 
    \item ${\displaystyle C(x) \wedge \neg E_C(x) \wedge \exists y
        \phi_C^\h(x,y) \rightarrow \forall z (\phi_C^\h(x,z) \rightarrow
        \g (x) = \g (z))}$.
    \end{enumerate}
  \item For each $C \in \phitrans$ and $(\g,\h) \in \{(\f,\b),
    (\b,\f)\}$ the formulas
    \begin{enumerate}
    \item ${\displaystyle E_{C^\g}\left(r^\g_C\right) \vee r^\g_C =
        \min(C^\g) \vee r^\g_C=\max(C^\g)}$,
    \item ${\displaystyle E_{C^\g}\left(s^\g_C\right) \vee s^\g_C =
        \min(C^\g) \vee s^\g_C=\max(C^\g)}$,
    \item ${\displaystyle r^\g_C \leq_{C^\g} s^\g_C}$,
    \item ${\displaystyle E_{C^\g}(x) \rightarrow \left(C(\h (x))
          \leftrightarrow r^\g_C <_{C^\g}x <_{C^\g} s^\g _C\right)}$.
    \end{enumerate}
  \item For each $m,n \in \NN$, $C \in \phiopen$, $D \in
    \phitrans$ and $\g \in \{\f,\b\}$ the formulas
    \begin{enumerate}
    \item ${\displaystyle E_C(x) \wedge E_C(\g^m(x)) \wedge \g^n(x)=x
        \rightarrow \g^m(x)=x}$,
    \item ${\displaystyle D(x) \wedge D(\g^m(x)) \wedge \g^n(x)=x \rightarrow
        \g^m(x) =x}$.
    \end{enumerate}
  \item For each $m \in \NN$ and $\g \in \{\f,\b\}$ the formula
    $\g^m(x)=x \rightarrow \g^{N_\Phi}(x)=x$.
  \end{lenumerate}
  This completes our list of axioms for $T(\Phi)$.
\end{df}

Our choice of axioms above and Sections \ref{foliation_orderings} and
\ref{progression_map} imply the following:

\begin{prop}
  \label{model1} 
  Let $\xi$ be a vector field on $\RR^2$ of class $C^1$ and definable
  in an \hbox{o-minimal} expansion of the real field, and let
  $\M_\xi$ be an $\la(\Phi_\xi)$-structure associated to
  $\xi$ as in Example \ref{omega_structure}.  Then $\M_{\xi}
  \models T(\Phi_{\xi})$. \qed
\end{prop}

\begin{df}
  \label{limit_cycle_set}
  We write $$\Phi_1 := \phitrans \cup \set{E_C:\ C \in \phiopen}.$$
  The following $\la(\Phi)$-formulas are of particular interest: for
  $C \in \Phi_1$, we let $\fix_C(x)$ be the formula $C(x) \wedge
  \f^{N_\Phi} (x) = x$ and $\fix_C(x,y)$ be the formula
  \begin{equation*}
    \exists z ((x\le_Cz\le_Cy \vee y\le_Cz\le_Cx) \wedge \fix_C(z)). 
  \end{equation*}
  Next, we let $\fixbd_C(x)$ be the formula
  \begin{equation*}
    \fix_C(x) \wedge \forall y \forall z \big(y<_Cx<_Cz \to \exists w
    (y<_Cw<_Cz \wedge \neg \fix_C(w))\big),
  \end{equation*}
  and let $\Lim_C(x)$ be the formula
  \begin{equation*}
    \fix_C(x) \wedge \exists y (C(y) \wedge y \ne x \wedge
    \neg\fix_C(x,y)). 
  \end{equation*}
\end{df}

\begin{expl}
  \label{omega_cycles}
  Let $\xi$ be a vector field on $\RR^2$ of class $C^1$ and definable
  in an \hbox{o-minimal} expansion of the real field, and let $\M_\xi$
  be an $\la(\Phi_\xi)$-structure associated to $\xi$ as in Example
  \ref{omega_structure}.  Let also $C \in \C_1:= \Ctrans \cup \{E_F:\
  F \in \Copen\}$.  Then the set $\fix_C(M)$ is the set of points in
  $C$ that belong to a cycle of $\xi$, the set $\fixbd_C(M)$ is the
  set of points in $C$ that belong to a boundary cycle of $\xi$, and
  the set $\Lim_C(M)$ is the set of points in $C$ that belong to a
  limit cycle of $\xi$.  Note that if $\xi$ is analytic, then the set
  $\fixbd_C(M)$ is discrete by Poincar\'e's Theorem \cite{poincare}
  (see also \cite[p. 217]{perko}); in particular, $\fixbd_C(M) =
  \Lim_C(M)$ in this case.

  In general, by Proposition \ref{fixed_points}, the cardinality of
  $\fixbd_C(M)$ is equal to the number of boundary cycles of $\xi$
  that intersect $C$.  Since every cycle of $\xi$ intersects the
  set $\bigcup \Ctan \cup \bigcup \Ctrans \cup \bigcup \Csingle$, it
  follows that, with $b(\xi)$ denoting the cardinality of the set
  of all boundary cycles of $\xi$, we have
  \begin{equation*}
    |\fixbd_C(M)| \le b(\xi) \le |\Ctan| + |\Csingle| + \sum_{D \in
      \Ctrans} |\fixbd_D(M)|.
  \end{equation*}
\end{expl}

\section{Iterating the progression maps}  \label{towards}

We continue to work with a flow configuration $\Phi$ as in Definition
\ref{config}.  Throughout this section, we fix $(\g,\h) \in
\{(\f,\b), (\b,\f)\}$.

For the next lemma, we denote by $\Theta_{(\g,\h)}$ the universal
closure of the conjunction of the formulas $(\bigwedge_{C \in \Phi_0}
\neg C(x)) \rightarrow \g(\h (x))=x$,
\begin{equation*}
  (C(x) \wedge E_C (\h (x))) \rightarrow \g(\h (x))=\g (x)
\end{equation*}
and
\begin{equation*}
  (E_C(x) \wedge \h (x) \ne s) \rightarrow \g(\h (x))=x
\end{equation*}
for each $C \in \phiopen$,
\begin{equation*}
  (C(x) \wedge \h (x) = e_C) \rightarrow \g(\h (x))=\g (x)
\end{equation*}
and 
\begin{equation*}
  (x=e_C \wedge \h (x) \ne s) \rightarrow \g(\h (x))=x
\end{equation*}
for each $C \in \phitan$, and $C(x) \rightarrow \g(\h (x))=x$ for each $C \in
\phitrans \cup \phisingle$.

\begin{lemma}
  \label{term1}
  $T(\Phi) \vdash \Theta_{(\g,\h)}$.
\end{lemma}

\begin{proof}
  Let $\M \models T(\Phi)$, and let $a \in M$ be such that $a \notin
  \bigcup_{C \in \Phi_0} C$.  Then by (F1), either $a = c$ for some $c
  \in \phisingle$, or $a = s$.  In the latter case, we have
  $\g(\h(a)) = \h(\g(a)) = a$ by (F5), so we may assume that
  $a = c$ for some $c \in \phisingle$.  Then $\h(\g(a)) =
  \g(\h(a)) = a$ by (F7)--(F9).

  The proofs of the other conjuncts is similar, using also (F12); we
  leave the details to the reader.
\end{proof}

\begin{cor}
  \label{term2} 
  Let $\phi$ be any quantifier-free $\la(\Phi)$-formula.  Then $\phi$
  is equivalent in $T(\Phi)$ to a quantifier-free formula $\phi'$ such
  that no term occurring in $\phi'$ contains both the symbols $\f$ and
  $\b$.
\end{cor}

\begin{proof}
  By induction on $l:= \max\{\length(t):\ t$ is a term occurring in
  $\phi\}$, using Lemma \ref{term1}.
\end{proof}

For the remainder of this section, we fix an arbitrary model $\M$ of
$T(\Phi)$.  To simplify notation, we omit the superscript $\M$ below
and write $\bar{C}:= C \cup \{\min(C),\max(C)\}$ for $C \in \Phi_1$.

\begin{df}
  \label{preimage}
  Let $C \in \Phi_1$ and $k \in \NN$.  We define
  \begin{equation*}
    G_C^k:= \{\g^l(z):\ z \text{ is a constant,}\ 0 \le
    l \le k \text{ and } \g^l(z) \in C\},
  \end{equation*}
  and we let $\O_C^k$ be the collection of all possible order types of
  pairs $(a,b) \in \bar{C}^{\,2}$ over $G_C^k$.  In addition, for
  $\zeta_0, \zeta_1 \in \bar{C}$ and $D \in \Phi_1$, we put
  \begin{equation*}
    \g^{-k}_D(\zeta_0,\zeta_1) := \set{x \in D:\ 
      \zeta_0 <_C \g^k(x) <_C \zeta_1}
  \end{equation*}
  and
  \begin{multline*}
    H^k_D(\zeta_0,\zeta_1):= \{\h^l(z):\ z \in
    \{\zeta_0, \zeta_1\} \text{ or } z \text{ is a constant}, \\
    0 \le l \le k \text{ and } \h^l(z) \in D\}.
  \end{multline*}
  Note that $G^k_C$ and $H^k_D(\zeta_0,\zeta_1)$, and hence $\O^k_C$,
  are finite sets whose cardinality is bounded by a number depending
  only on the language and $k$, but independent of $\M$, $C$, $D$,
  $\zeta_0$ or $\zeta_1$.
\end{df}

\begin{prop}
  \label{pb2} 
  Let $C, D \in \Phi_1$, $\zeta_0,\zeta_1 \in \bar{C}$ and $k \in
  \NN$.
  \begin{enumerate}
  \item The set $\g^{-k}_D(\zeta_0,\zeta_1)$ is a union of points in
    $H^k_D(\zeta_0, \zeta_1)$ and open intervals with endpoints in
    $H^k_D(\zeta_0,\zeta_1)$.
  \item For each $\vartheta \in \O^k_C$, there is a conjunction
    $\sigma_\vartheta(x,y_0,y_1)$ of atomic formulas with free variables
    $x$, $y_0$ and $y_1$ such that whenever $(\zeta_0,\zeta_1)$ have
    order type $\vartheta$ over $G^k_C$, the set
    $\g^{-k}_D(\zeta_0,\zeta_1)$ is defined by the formula
    $\sigma_\vartheta(x,\zeta_0,\zeta_1)$.
  \item $\g^k$ restricted to $\g^{-k}_D(\zeta_0,\zeta_1)$ is
    continuous.
  \end{enumerate}
\end{prop}

\begin{proof}
  For every $x \in \g^{-k}_D(\zeta_0, \zeta_1)$, there is a
  sequence $E = (E_0, \dots, E_k)$ of elements of $\Phi_2:= \Phi_1
  \cup \set{\{c\}:\ c \in \phisingle} \cup \set{\{e_C\}:\ C \in
    \phitan}$ such that $E_0 = D$, $E_k = C$ and $\g^i(x) \in E_i$
  for $i=0, \dots, k$.  Thus, we fix a sequence $E = (E_0, \dots, E_k)
  \in \Phi_2^{k+1}$ with $E_k = C$, and we define the set
  \begin{equation*}
    \g^{-k}_E(\zeta_0,\zeta_1) := \set{x \in M:\ \g^i(x) \in E_i \text{
        for } i=0, \dots, k,\ \zeta_0 <_C
      \g^k(x) <_C \zeta_1};
  \end{equation*}
  it suffices to prove the proposition with
  $\g^{-k}_E(\zeta_0,\zeta_1)$  and $H^k_{E_0}(\zeta_0,\zeta_1)$ in
  place of $\g^{-k}_D(\zeta_0,\zeta_1)$ and $H^k_D(\zeta_0,\zeta_1)$.

  Next, we note that if $E_i \in \set{\{c\}:\ c \in \phisingle} \cup
  \set{\{e_C\}:\ C \in \phitan}$ for some $i \in \{1, \dots, k-1\}$,
  then $a \in \g^{-k}_E(\zeta_0,\zeta_1)$ if and only if $\g^i(a)$ is
  the unique constant in $E_i$ and $\zeta_0 <_C \g^k(a) <_C \zeta_1$,
  so the proposition follows in this case.  

  We therefore assume from now on that $E_i \in \Phi_1$ for each $i=0,
  \dots, k$, and in this case we prove the proposition with part (1)
  replaced by
  \begin{itemize}
  \item[(1)'] The set $\g^{-k}_E(\zeta_0,\zeta_1)$ is an open interval
    with endpoints in $H^k_{E_0}(\zeta_0,\zeta_1)$.
  \end{itemize}
  We proceed by induction on $k$.  The case $k=0$ is trivial, so we
  assume that $k>1$.  By Axiom (F8), the set
  $\g^{-1}_{(E_{k-1},E_k)}(\zeta_0,\zeta_1)$ is an open interval whose
  endpoints $\eta_0, \eta_1$ belong to the set
  $H^1_{E_{k-1}}(\zeta_0,\zeta_1)$ and are determined by the order
  type of $(\zeta_0,\zeta_1)$ over $G^1_{E_k}$.  In fact, we claim
  that the order type of $(\eta_0,\eta_1)$ over $G^{k-1}_{E_{k-1}}$ is
  determined by the order type of $(\zeta_0,\zeta_1)$ over
  $G^k_{E_k}$; together with the inductive hypothesis applied to
  $\g^{k-1}_{(E_0,\dots,E_{k-1})}(\eta_0,\eta_1)$, the proposition then
  follows, because $H^{k-1}_{E_0}(c,d)$ is contained in
  $H^k_{E_0}(\zeta_0,\zeta_1)$ for all $c,d \in
  H^1_{E_{k-1}}(\zeta_0,\zeta_1)$.

  To see the claim, assume first that $E_k = E_C$ for some $C \in
  \phiopen$.  Then by Axiom (F8), the set $\{\g(z):\ z \in
  G^{k-1}_{E_{k-1}}\}$ is contained in $G^k_{E_k}$ and the claim
  follows in this case.  So we assume that $E_k \in \phitrans$.  Then
  by Axiom (F13), $E_{k-1} = E_C$ for some $C \in \phiopen$ and there
  are constants $a$ and $b$ such that
  \begin{equation*}
    (\eta_0,\eta_1) \subseteq (a,b) = \g^{-1}(E_k) = h(E_k)
    \quad\text{(as intervals)}.
  \end{equation*}
  Hence the order type of $(\eta_0,\eta_1)$ over $G^{k-1}_{E_C}$ is
  determined by the order type of $(\eta_0,\eta_1)$ over the set $G'
  := \set{z \in G^{k-1}_{E_C}:\ a <_C z <_C b}.$ Then again by Axiom
  (F8), the set $\{\g(z):\ z \in G'\}$ is contained in $G^k_{E_k}$ and
  the claim also follows in this case.
\end{proof}

\begin{cor}
  \label{def_complete}
  Let $C \in \Phi_1$ and put $G := \g^{-N}_C(\min(C),\max(C))$.
  \begin{enumerate}
  \item The set $\fixbd_C(M)$ is a closed and nowhere dense subset of
    $G$. 
  \item Assume that $\Phi = \Phi_\xi$ and $\M \equiv \M_\xi$ for
    some definable vector field $\xi$ of class $C^1$ on $\RR^2$.  Then
    for every $c \in G \setminus \fixbd_C(M)$, there are $a, b \in
    \bar{C}$ such that
    \begin{equation*}
      a = \sup\set{x \in \fixbd_C(M) \cup (\bar{C} \setminus
        G):\ x <_C c}
    \end{equation*}
    and
    \begin{equation*}
      b = \inf\set{x \in \fixbd_C(M) \cup (\bar{C} \setminus G):\ c
        <_C x}. 
    \end{equation*}
  \end{enumerate}
\end{cor}

\begin{proof}
  Part (1) follows from the continuity of $\g^N\rest{G}$ and the
  definition of the set $\fixbd_C(M)$.  Part (2) follows from part (1)
  and the fact that $C^{\M_\xi}$ is complete.
\end{proof}

Finally, for each $C \in \Phi_1$ we let $\bar{C}(x)$ abbreviate $C(x) \vee
x=\min(C) \vee x=\max(C)$.  We let $G^k$ be the set of all
$\la(\Phi)$-terms $\g^jc$ such that $0 \le j \le k$ and $c$ is a
constant symbol, and we let $\O^k$ be the set of all formulas of the
form
\begin{equation*}
  \left(\bar{C}(y_0) \wedge \bar{C}(y_1)\right) \wedge
  \bigwedge_{\{\tau,\rho\} \subseteq G^k 
    \cup \{y_0,y_1\}} (\tau \ast_{\{\tau,\rho\}} \rho), 
\end{equation*}
where $C \in \Phi_1$ and $\ast_{\{\tau,\rho\}} \in \{<_C, >_C, =,
\ne\}$.  The cardinalities of $G^k$ and $\O^k$ are bounded by a number
depending only on $k$ (and on $\la(\Phi)$).  Moreover in $\M$, each
formula $\vartheta \in \O^k$ determines an order type in $\O^k_C$, for
some $C \in \Phi_1$; and conversely, every order type in $\O^k_C$ with
$C \in \Phi_1$ is determined by some formula $\vartheta \in \O^k$.
Thus we obtain the following from Proposition \ref{pb2}:

\begin{cor}
  \label{shape} 
  Let $k \in \NN$.  Then there are $l = l(k) \in \NN$ and
  quantifier-free formulas $\vartheta^k_1(y_0,y_1), \dots,
  \vartheta^k_l(y_0,y_1)$ with free variables $y_0$ and $y_1$ such that
  \begin{enumerate}
  \item ${\displaystyle T(\Phi) \vdash \bigvee_{i=1}^l
      \vartheta^k_i(y_0,y_1) \leftrightarrow \bigvee_{C \in
        \Phi_1}\left(\bar{C}(y_0) \wedge \bar{C}(y_1)\right)}$;
  \item for every $D \in \Phi_1$ there are quantifier-free formulas
    $\sigma^{D,k}_i(x,y_0,y_1)$ with free variables $x$, $y_0$ and
    $y_1$, $i = 1, \dots, l$, such that if $\M \models
    \vartheta^k_i(\zeta_0,\zeta_1)$ for $\zeta_0, \zeta_1 \in M$ and
    some $i$, then the set $\g^{-k}_D(\zeta_0,\zeta_1)$ is defined by
    the formula $\sigma^{D,k}_i(x,\zeta_0,\zeta_1)$. \qed
  \end{enumerate}
\end{cor}

\begin{nrmk} 
  \label{shape_remark}
  We obtain analogous statements to Proposition \ref{pb2} and
  Corollary \ref{shape} if we replace the open interval
  $(\zeta_0,\zeta_1)$ by a half-open or closed interval.
\end{nrmk}

\section{Dulac flow configurations}  \label{dulac_flow}

It is clear from Remark \ref{omega_cycles} that, for a vector field
$\xi$ on $\RR^2$ definable in $\R$, the set of boundary cycles of
$\xi$ is represented in $\M_\xi$ by the definable sets
$\fixbd_C(M)$.  The following example shows that the theory $T(\Phi)$
has hardly any implications for the nature of these sets.

\begin{expl}
  \label{perturb}
  Consider the vector field $\zeta$ of Example \ref{radial}, and let
  $\C$ be the piecewise trivial decomposition obtained there.  We
  denote by $\Phi_\zeta$ the flow configuration corresponding to this
  $\C$ and write $$C_0:= \set{(x,y):\ x>0,\ y=0} \in \C.$$  We show
  here how to define, given any closed and nowhere dense subset $F$ of
  $C_0$, a vector field $\zeta'$ of class $C^\infty$ for which
  $\Phi_\zeta$ is still a flow configuration and such that
  $\fixbd_{C_0}(M_{\zeta'}) = F$.

  First, given $0<a<b<\infty$, we let $d_{(a,b)}: \RR^2 \into \RR$ be
  the function $d_{(a,b)}(x,y) := (b^2 - (x^2+y^2))((x^2+y^2) - a^2)$,
  and we let $e_{(a,b)}: \RR^2 \into \RR$ be the $C^\infty$ function
  defined by $e_{(a,b)}(x,y) := \exp(-1/d_{(a,b)}(x,y))$.  We let
  $\zeta_{(a,b)}$ be the vector field of class $C^\infty$ on the annulus
  $A_{(a,b)}:= \set{(x,y): d_{(a,b)}(x,y) > 0}$ defined by
    \begin{equation*}
      \zeta_{(a,b)} := -\left(y + e_{(a,b)}(x,y) x\right)
      \frac{\partial}{\partial x}+ \left(x - e_{(a,b)}(x,y) y\right)
      \frac{\partial}{\partial y}. 
    \end{equation*}

  Second, let $F \subseteq C_0$ be an arbitrary closed and nowhere
  dense subset.  Then $C_0 \setminus F$ is open in $C_0$ and hence the
  union of countably many disjoint open intervals $I_0, I_1, I_2,
  \dots$.  We let $\zeta'$ be the vector field on $\RR^2$ of class
  $C^\infty$ defined by
  \begin{equation*}
    \zeta'(x,y) :=
    \begin{cases}
      \zeta_{I_j}(x,y) &\text{if } (x,y) \in A_{I_j} \text{ for some } j
      \in \NN, \\ \zeta(x,y) &\text{otherwise.}
    \end{cases}
  \end{equation*}
  (Note that by Wilkie's Theorem \cite{wilkie}, $\zeta'$ is
  definable in some o-minimal expansion of the real field if and only
  if $F$ is finite.)  
\end{expl}

In view of the previous example, we now introduce a strengthening of
the setting described in Section \ref{flowconfig}.

\begin{df} 
  \label{dulac_configuration}
  A \textbf{Dulac flow configuration} $\Psi$ is a pair $(\Phi, \nu)$
  such that $\Phi$ is a flow configuration and $\nu \in \NN$.
\end{df} 

\begin{expl}
  \label{dulac_omega_config}
  Let $\xi$ be a definable vector field on $\RR^2$ of class $C^1$.
  Let $\Phi=\Phi_{\xi}$ be a flow configuration associated to $\xi$ as
  in Example \ref{omega_config} and let $\M_{\xi}$ be the associated
  $\la(\Phi_\xi)$-structure described in Example
  \ref{omega_structure}.  Assume that there is a $\nu \in \NN$ such
  that for each $C \in \Phi_1$, the set $\fixbd_C(M_{\xi})$ has
  cardinality at most $\nu$.  Then $\Psi_{\xi} := (\Phi_{\xi},\nu)$ is
  called a \textbf{Dulac flow configuration associated to $\xi$}.
\end{expl}

For the remainder of this section, we fix a Dulac flow configuration
$\Psi = (\Phi,\nu)$.

\begin{df}
  \label{dulac_language}
  The language $\la(\Psi)$ consists of the symbols of $\la(\Phi)$
  together with the following symbols for each $C \in \Phi_1$:
  \begin{renumerate}
  \item binary predicates $R_C$ and $S^\f_{m,C}$, $B^\f_{m,C}$,
    $S^\b_{m,C}$ and $B^\b_{m,C}$ for each $m \in \NN$;
  \item constant symbols $\gamma^1_C, \dots, \gamma^\nu_C$.
  \end{renumerate}
  We put $\Gamma = \Gamma(\Psi) := \set{\gamma^j_C:\ C \in \Phi_1,\
    j=1, \dots, \nu}$.
\end{df}

\begin{expl}
  \label{dulac_omega_structure}
  Let $\xi$ be a definable vector field on $\RR^2$ of class $C^1$, and
  let $\M_\xi$ be an $\la(\Phi_\xi)$-structure associated to
  $\xi$ as in Example \ref{omega_structure}.  Assume that there is
  a $\nu \in \NN$ such that for each $C \in \Ctrans \cup \Copen$,
  the set $\fixbd_C(M_{\xi})$ has cardinality at most $\nu$, and
  let $\Psi_\xi$ be a Dulac flow configuration associated to
  $\xi$ as in Example \ref{dulac_omega_config}.  We expand
  $\M_\xi$ into an $\la(\Psi_\xi)$-structure $\M^D_\xi$ as
  follows: for each $C \in \Phi_1$,
  \begin{renumerate}
  \item $R_C$ is interpreted as the set
    \begin{equation*}
      \set{(x,y) \in \bar{C}^{\,2}:\ \exists z
        (x<_Cz<_Cy \wedge \fix_C(z)) \vee \left(x=y \wedge 
          \fix_C(x)\right)};
    \end{equation*}
  \item for $m \in \NN$, $\g \in \{\f,\b\}$ and $G \in \{S^\g_{m,C},
    B^\g_{m,C}\}$, we put
    \begin{equation*}
      \ast :=
      \begin{cases}
        <_C &\text{if } G \text{ is } S^\g_{m,C}, \\ >_C
        &\text{if } G \text{ is } B^\g_{m,C},
      \end{cases}
    \end{equation*}
    and we interpret $G$ as the union of the sets
    \begin{equation*}
      \set{(x,y) \in \bar{C}^{\,2}:\ \exists z \big(C(z) \wedge x <_C z <_C
        y \wedge C (\g^m(z)) \wedge \g^m(z) \ast z\big)}
    \end{equation*}
    and the set $\set{(x,x): C(x) \wedge C(\g^m(x)) \wedge \g^m(x) \ast 
      x}$;
  \item if $a_1<_C \cdots <_Ca_m$ are the points in $C$ that lie on
    boundary cycles of $\xi$, we interpret $\gamma^j_C$ as $a_j$ if
    $1 \le j \le m$ and as $\max(C)$ if $m < j \le \nu$.
  \end{renumerate}
  This completes the description of $\M^D_{\xi}$.
\end{expl}

\begin{df}
  \label{dulac_theory}
  Inspired by the previous example, we let $T(\Psi)$ be the
  $\la(\Psi)$-theory consisting of $T(\Phi)$ and the universal
  closures of the formulas in the axiom schemes (D1)--(D6) below.
  \begin{lenumerate}{D}{0}
  \item For each $C \in \Phi_1$, $m \in \NN$ and $G \in
    \{R_C,S^\f_{m,C}, B^\f_{m,C}, S^\b_{m,C}, B^\b_{m,C}\}$, the
    formulas
    \begin{enumerate}
    \item ${\displaystyle G(x,y) \rightarrow \left(\bar{C}(x) \wedge
          \bar{C}(y)\right)}$,
    \item ${\displaystyle G(x,y) \rightarrow \left(x \leq_C y \vee
          (x=\min(C) \wedge y=\max(C))\right)}$.
    \end{enumerate}
  \item For each $C \in \Phi_1$ the formulas
    \begin{enumerate}
    \item $R_C(x,y) \leftrightarrow \exists z (x <_C z <_C y \wedge
      \fix_C(z))$, and
    \item $R_C(x,x) \leftrightarrow \fix_C(x)$.
    \end{enumerate}
  \item For each $m \in \NN$, $C \in \Phi_1$ and $\g \in \{\f,\b\}$
    the formulas
    \begin{enumerate}
    \item $S^\g_{m,C}(x,y) \leftrightarrow \exists z (x <_C z <_C y \wedge
      \g^m(z) <_C z)$,
    \item $S^\g_{m,C}(x,x) \leftrightarrow (C(x) \wedge \g^m (x) <_C x)$,
    \item $B^\g_{m,C}(x,y) \leftrightarrow \exists z (x <_C z <_C y
      \wedge z <_C \g^m(z))$,
    \item $B^\g_{m,C}(x,x) \leftrightarrow (C(x) \wedge x <_C \g^m (x))$.
    \end{enumerate}
  \item For each $m \in \NN$, $C \in \Phi_1$, $\g \in \{\f,\b\}$ and
    $G \in \{R_C, B^\g_{m,C}, S^\g_{m,C}\}$ the formula
    \begin{multline*}
      \Big[(G(x,y) \wedge \forall z \left(x <_C z <_C y \to \bar{C} (\g^m
          (z))\right) \\ \wedge \neg\exists z \left(x <_C z <_C y \wedge
          \fixbd_C(z)\right)\Big] \\ \to \forall z (x <_C z <_C y
      \to G(z,z)).
    \end{multline*}
  \item[(D5)$_\nu$] For each $C \in \Phi_1$ the formulas
    \begin{enumerate}
    \item ${\displaystyle \bar{C}\left(\gamma_C^j\right) \wedge
        \left(C\left(\gamma^j_C\right) \rightarrow
          \fix_C\left(\gamma_C^j\right)\right)}$ for $j=0, \dots, \nu$,
    \item ${\displaystyle \gamma_C^j \leq_C \gamma^{j+1}_C} \wedge
      \left(\gamma_C^j = \gamma_C^{j+1} \to \gamma_C^j =
        \max(C)\right)$ for $j=0, \dots, \nu-1$.
    \end{enumerate}
  \item[(D6)$_\nu$] For each $C \in \Phi_1$ the formula
    \begin{equation*}
      (C(x) \wedge \fixbd_C(x)) \leftrightarrow
      \bigvee_{j=1}^\nu \left(x=\gamma^j_C \wedge
        C\left(\gamma^j_C\right)\right).
    \end{equation*}
  \end{lenumerate}
  This completes the description of the axioms.
\end{df}

\begin{prop}  
  \label{dulacprop}
  If $\xi$ is a definable vector field on $\RR^2$ of class $C^1$ with
  finitely many boundary cycles, then $\M_{\xi}^D \models
  T(\Psi_\xi)$.
\end{prop}

\begin{proof}  
  This is almost immediate from the definition of $\M^D_{\xi}$ and
  Proposition \ref{model1}, except perhaps for Axiom (D4), which
  follows from Proposition \ref{pb2} and the fact that every bounded
  subset of $\RR$ has an infimum.
\end{proof} 

\begin{nrmk}
  \label{predulac}
  Let $T(\Phi)'$ be the union of $T(\Phi)$ with Axioms (D1)--(D4)
  only.  Since (D1)--(D3) just extend $T(\Phi)$ by definitions in the
  sense of Section 4.6 in Shoenfield \cite{shoenfield}, the argument
  in the proof of the previous proposition shows that any
  $\la(\Phi_\xi)$-structure $\M_\xi$ as defined in Example
  \ref{omega_structure} can be expanded to a model $\M'_\xi$ of
  $T(\Phi)'$.
\end{nrmk}

\section{Quantifier elimination for $T(\Psi)$}  \label{QE}

We fix a Dulac flow configuration $\Psi = (\Phi,\nu)$; our ultimate
goal is to show that $T(\Psi)$ eliminates quantifiers.  Most of the
work in this section goes towards showing that, in order to eliminate
quantifiers, we need only consider formulas of the form $\exists y
\phi(x,y)$ where $\phi$ is of a special form.

\begin{term}
  Let $x = (x_1, \dots, x_m)$ be a tuple of variables and $y$ and $z$
  single variables.  To simplify terminology, we write ``term'' and
  ``formula'' for ``$\la(\Psi)$-term'' and ``$\la(\Psi)$-formula''.
  For a formula $\phi$, we write $\phi(x,y)$ to indicate that the free
  variables of $\phi$ are among $x_1, \dots, x_m$ and $y$.  A
  \textbf{binary atomic formula} is a formula of the form $A t_1 t_2$,
  where $A$ is a binary relation symbol in $\la(\Psi)$ and $t_1$ and
  $t_2$ are terms.
\end{term}

For this section fix an arbitrary model $\M$ of $T(\Psi)$; again, we
omit the superscript $\M$ when interpreting predicates in $\M$.

\begin{df} 
  \label{orderformula}
  An \textbf{order formula} is a quantifier-free $\la(\Phi) \cup
  \Gamma$-formula.  A \textbf{$z$-order formula} is a quantifier-free
  formula $\phi$ such that every atomic subformula of $\phi$
  containing $z$ is an $\la(\Phi) \cup \Gamma$-formula.

  A $z$-order formula $\phi$ is \textbf{minimal} if the only subterm
  of $\phi$ containing $z$ is $z$ itself and every binary atomic
  subformula $At_1t_2$ of $\phi$ is such that at most one of $t_1$ and
  $t_2$ contains $z$.
\end{df}

Our first goal is to show that we may, in order to prove quantifier
elimination, restrict our attention to $y$-order formulas.  This
argument is based on the following lemma, which will also be of use
later. 

\begin{lemma}
  \label{rsimp}
  Let $G \in \la(\Psi) \setminus \la(\Phi)$.
  \begin{enumerate}
  \item The formula $Gyy$ is equivalent in $T(\Psi)$ to a minimal
    $y$-order formula $\psi(y)$.
  \item The formula $Gyz$ is equivalent in $T(\Psi)$ to a formula
    $\psi(y,z)$ that is both a minimal $y$-order formula and a
    minimal $z$-order formula.
  \end{enumerate}
\end{lemma}

\begin{proof}
  Let $C \in \Phi_1$, $m \in \NN$ and $\g \in \{\f,\b\}$ be such that
  $G$ is one of $R_C$, $S^\g_{m,C}$ or $B^\g_{m,C}$.  In this proof, we write
  $<$ instead of $<_C$; if $G$ is $R_C$, we assume $m=N = N_\Phi$.  By
  Corollary \ref{shape}(1), any formula $\phi$ is equivalent in
  $T(\Psi)$ to the conjunction of the formulas $\vartheta_i
  \rightarrow \phi$, where $i \in \{1, \dots, l(m)\}$ and
  $\vartheta_i$ is the formula $\vartheta^m_i(\min(C),\max(C))$.
  Hence it suffices to prove the lemma with each $\vartheta_i
  \rightarrow G(y,y)$ in place of $G(y,y)$ and each $\vartheta_i \to
  G(y,z)$ in place of $G(y,z)$; so we also fix an $i$ below and write
  $\vartheta$ in place of $\vartheta_i$.  Now by Corollary
  \ref{shape}(2), there are finitely many terms
  $\alpha_j^0,\alpha_j^1$ for $1 \leq j \leq r$, built up exclusively
  from constants, such that whenever $\M \models \vartheta$ the set
  $\{z \in C : \g^m(z) \in C\}$ is the union of the open intervals
  $I_j = (\alpha_j^0,\alpha_j^1)$ and points $\alpha_j^0=\alpha_j^1$.

  (1) We claim that the formula $\vartheta \rightarrow G(y,y)$ is
  equivalent to $\vartheta \to \psi^G$, where $\psi^G$ is of the form
  \begin{equation*}
    C(y) \wedge \left(\bigvee_{1 \leq j \leq r} (\alpha^0_j < y <
      \alpha_j^1 \vee \alpha_j^0=y=\alpha_j^1)\right) \wedge
    \left( \bigvee_{\beta \in Y} \psi_\beta^G
      \vee \bigvee_{\beta_0,\beta_1 \in Y} \psi_{\beta_0,\beta_1}^G  \right)
  \end{equation*}
  with $Y := \Gamma \cup \{\alpha_j^l : l \in \{0,1\} \mbox{ and } 1
  \leq j \leq r\}$, and for each $\beta \in Y$, the formula
  $\psi^G_\beta$ is $C(y) \wedge ((y=\beta \wedge G(\beta,\beta)) \vee y =
  t^G)$ with 
  \begin{equation*}
    t^G \text{ the term } 
    \begin{cases}
      y &\text{if } G \text{ is } R_C, \\ \h^{m}\min(C) &\text{if } G
      \text{ is } S^\g_{m,C}, \\ \h^{m}\max(C) &\text{if } G \text{ is
      } B^\g_{m,C},
    \end{cases}
  \end{equation*}
  and for each $\beta_0,\beta_1 \in Y$, the formula
  $\psi^G_{\beta_0,\beta_1}$ is of the form
  \begin{equation*}
    (C(\beta_0) \vee
    \beta_0=\min(C)) \wedge (C(\beta_1) \vee \beta_1=\max(C)) \wedge
    \beta_0 < y < \beta_1 \wedge \eta^G_{\beta_0,\beta_1},
  \end{equation*}
  where
    \begin{equation*}
     \eta^G_{\beta_0,\beta_1} \text{ is } 
      \begin{cases}
        \neg S^\g_{N,C} (\beta_0,\beta_1) \wedge \neg B^\g_{N,C}
        (\beta_0,\beta_1) &\text{if } G \text{ is } R_C, \\ \neg
        B^\g_{m,C} (\beta_0,\beta_1) \wedge \neg R_C (\beta_0,\beta_1)
        &\text{if } G \text{ is } S^\g_{m,C}, \\ \neg S^\g_{m,C}
        (\beta_0,\beta_1) \wedge \neg R_C (\beta_0,\beta_1) &\text{if } G
        \text{ is } B^\g_{m,C}.
      \end{cases}
    \end{equation*}
  Note that $\vartheta \to \psi^G$ is a minimal $y$-order formula;
  thus, the proof of part (1) is finished once we prove the claim.

  We prove the claim for $R_C$; the other cases of $G$ are similar and
  left to the reader.  Suppose that $\M \models \vartheta$ and pick an
  $a \in M$ such that $\M \models R_C(a,a)$.  Then $\M \models
  \alpha_j^0 \leq a \leq \alpha_j^1$ for some $j \in \{1, \dots r\}$.
  If $a=\beta$ for some $\beta \in Y$, we are done, so we assume $a
  \ne \beta$ for all $\beta \in Y$.  Then there are $\beta_0, \beta_1
  \in Y$ such that $\M \models \beta_0 < a < \beta_1$ and $\M \models
  \neg(\beta_0 < \beta < \beta_1)$ for every $\beta \in Y$.  Hence by
  Axiom (D4), $\M \models R_C(b,b)$ for every $b \in
  (\beta_0,\beta_1)$, so $\M \models \neg S^\g_{m,C}(\beta_0,\beta_1)
  \wedge \neg B^\g_{m,C}(\beta_0,\beta_1)$ as required.  The converse
  of the claim is immediate.

  (2) The formula $\vartheta \to G(y,z)$ is in turn equivalent in
  $T(\Psi)$ to 
  \begin{equation*}
    \vartheta \to (G(y,z) \wedge (y=\min(C) \vee y = \max(C) \vee
    C(y))); 
  \end{equation*}
  since the lemma is immediate for the formulas $\vartheta \to (G(y,z)
  \wedge y=\min(C))$ and $\vartheta \to (Gyz \wedge y = \max(C))$, we
  need only consider $\vartheta \to (G(y,z) \wedge C(y))$.  We claim
  that the latter is equivalent to $\vartheta \to \psi^G$, where
  $\psi^G$ is of the form
  \begin{equation*}
    C(y) \wedge (C(z) \vee z=\max(C)) \wedge y \le z \wedge \left((y=z
    \wedge G(y,y)) \vee \left(y<z \wedge \eta^G\right)\right),
  \end{equation*}
  $\eta^G$ is the formula
  \begin{equation*}
    \bigvee_{\beta \in Y}
    (y=\beta \wedge G (\beta, z)) \vee \bigvee_{\beta \in Y}(y < \beta
    < z \wedge G(\beta,\beta)) \vee \bigvee_{\beta_o, \beta_1 \in
      Y,\ 1 \leq j \leq r} \eta^G_{\beta_0,\beta_1,j}
  \end{equation*}
  and for each $\beta_0,\beta_1 \in Y$ and $j \in \{1, \dots, r\}$,
  the formula $\eta^G_{\beta_0,\beta_1,j}$ is
  \begin{equation*}
    \beta_0 < y \wedge z < \beta_1 \wedge
    \alpha^0_j \leq \beta_0 \wedge \beta_1 \leq \alpha_j^1 \wedge
    G(\beta_0,\beta_1) \wedge \eta^G_{\beta_0,\beta_1}
  \end{equation*}
  with $\eta^G_{\beta_0,\beta_1}$ defined as for part (1). 

  We again prove the claim for $R_C$, leaving the other cases of $G$
  to the reader.  Suppose that $\M \models \vartheta$ and $\M \models
  R_C(a,b) \wedge C(b)$ and work inside $\M$.  Suppose that $a \ne
  \beta$ for all $\beta \in Y$ and that $\M \models \neg (a < \beta <
  b \wedge R_C(\beta,\beta))$ for every $\beta \in Y$.  Then $\f^N(d)
  = d$ for some $d \in (a,b)$, and $d \in (\alpha^0_j,\alpha^1_j)$ for
  some $j$.  Moreover, there are $\beta_0, \beta_1 \in Y$ such that $d
  \in (\beta_0, \beta_1)$ and $\beta \notin (\beta_0, \beta_1)$ for
  every $\beta \in Y$.  Hence by Axiom (D4), we get $\M \models \neg
  S^\g_{N,C}(\beta_0,\beta_1) \wedge \neg
  B^\g_{N,C}(\beta_0,\beta_1)$, as required.  The converse of the
  claim is straightforward.

  By symmetry, a similar claim holds with $\vartheta \to (G(y,z)
  \wedge C(z))$ in place of $\vartheta \to (G(y,z) \wedge C(y))$.
  Combining these two claims with part (1) now yields part(2).
\end{proof}

\begin{cor}
  \label{fsimp} 
  Every quantifier-free formula $\phi(x,y)$ is equivalent in $T(\Psi)$
  to a $y$-order formula $\psi(x,y)$.
\end{cor}

\begin{proof}
  It suffices to prove the proposition for all atomic formulas; the
  relevant atomic formulas are handled in Lemma \ref{rsimp}.
\end{proof}

Our second goal of this section is to show that we only need consider,
for quantifier elimination, $y$-order formulas in which the complexity
of any term involving $y$ is as low as possible.  Minimal $y$-order
formulas are examples of such $y$-order formulas; but we cannot always
reduce to minimal $y$-order formulas.

\begin{df}
  \label{height}
  Let $t$ be a term.  The \textbf{$z$-height $h_z(t)$ of $t$} is defined
  as follows:
  \begin{renumerate}
    \item if $z$ does not occur in $t$, then $h_z(t):= 0$;
    \item $h_z(z) := 1$;
    \item if $t$ is $\f t'$ or $\b t'$ for some term $t'$ and $z$
      occurs in $t'$, then $h_z(t) := h_z(t') + 1$.
  \end{renumerate}

  Let $A(t_1,t_2)$ be a binary atomic formula; the \textbf{$z$-height
    $h_z(A(t_1,t_2))$ of $A(t_1,t_2)$} is defined as the pair
  $(a,b) \in \NN^2$, where
  \begin{equation*}
    a := 
    \begin{cases}
      1 &\text{if } z \text{ occurs in both } t_1 \text{ and } t_2, \\
      0 &\text{otherwise,}
    \end{cases}
  \end{equation*}
  and
  \begin{equation*}
    b :=
    \begin{cases}
      \min\{h_z(t_1),h_z(t_2)\} &\text{if } z \text{ occurs in both }
      t_1 \text{ and } t_2, \\ \max\{h_z(t_1),h_z(t_2)\}
      &\text{otherwise.} 
    \end{cases}
  \end{equation*}

  Let $B(t)$ be a unary atomic formula; the \textbf{$z$-height $h_z(B(t))$
    of $B(t)$} is defined by $h_z(B(t)):= (0,h_z(t)) \in \NN^2$.

  Let $\phi$ be a quantifier-free formula; the \textbf{$z$-height
    $h_z(\phi)$ of $\phi$} is the maximum of the set $\set{h_z(\psi):\
    \psi \text{ is an atomic subformula of } \phi}$ with respect to
  the lexicographic ordering of $\NN^2$.  We write $h_z(\phi) =
  (h_z^1(\phi),h_z^2(\phi))$ below.

  Finally, a term $t$ is \textbf{mixed} if it contains both function
  symbols $\f$ and $\b$; otherwise $t$ is called \textbf{unmixed}.
\end{df}

\begin{expl}
  \label{height-example}
  Let $\phi$ be a $z$-order formula.  Then $h_z(\phi) \le (0,1)$ if
  and only if $\phi$ is minimal.
\end{expl}

\begin{lemma}
  \label{unmix}
  Let $\phi(x,y)$ be a $y$-order formula.  Then there is a $y$-order
  formula $\psi(x,y)$ that contains no mixed terms such that $\phi$
  and $\psi$ are equivalent in $T(\Psi)$ and $h_y(\psi) \le h_y(\phi)$.
\end{lemma}

\begin{proof}
  Let $\phi'$ be the $\la(\Phi)$-formula obtained from $\phi$ by
  replacing each constant $\gamma^j_C$ by a new variable $z^j_C$, for
  $C \in \Phi_1$ and $j=1, \dots, \nu$.  By Lemma \ref{term1}, $\phi'$
  is equivalent in $T(\Phi)$ to a quantifier-free $\la(\Phi)$-formula
  $\psi'$ that is a disjunction of formulas of the form $\eta \wedge
  \xi$, where $\xi$ is obtained from $\phi'$ by replacing each mixed
  subterm by an unmixed term of lower $y$-height, and where $\eta$ is
  a conjunction of some of the premises of the implications occurring
  in $\Theta_{(\f,\b)}$ and in $\Theta_{(\b,\f)}$ with $x$ there
  replaced by various unmixed subterms of $\phi'$.  Clearly $h_y(\xi)
  \le h_y(\phi')$ for every such $\xi$; since $h_y^1(\eta) = 0$ for
  every such $\eta$, it follows that $h_y(\psi') \le h_y(\phi')$ if
  $h_y^1(\phi') = 1$.  On the other hand, if $h_y^1(\phi') = 0$, then
  every subterm $t$ of $\phi'$ satisfies $h_y(t) \le h_y^2(\phi')$; so
  $h_y(\eta) \le h_y(\phi')$ for every such $\eta$.  Therefore, we
  always have $h_y(\psi') \le h_y(\phi') = h_y(\phi)$, and we let
  $\psi$ be the $y$-order formula obtained from $\psi'$ by replacing
  each variable $z^j_C$ again by $\gamma^j_C$.
\end{proof}

Below we let $\iota(y)$ denote the formula $\bigwedge_{C \in \phiopen}
C(y) \to E_C(y)$ and we put $$T' := T(\Psi) \cup \{\iota(y)\}.$$

\begin{lemma}
  \label{height_simplification}
  Let $\phi(x,y)$ be a $y$-order formula.  Then there is a $y$-order
  formula $\psi(x,y)$ such that $\phi$ is equivalent in $T'$ to $\psi$
  and $h^2_y(\psi) \le 1$.
\end{lemma}

\begin{proof}
  By induction on $h_y(\phi)$; the case where $h_y^2(\phi) \le 1$ is
  trivial, so we assume that $h_y^2(\phi) > 1$ and we prove that
  \begin{itemize}
  \item [$(\ast)$] there exists an order formula $\psi(x,y)$ such that
    $\phi$ is equivalent in $T'$ to $\psi$ and $h_y(\psi) <
    h_y(\phi)$.
  \end{itemize}
  To do so, we fix arbitrary $(\g,\h) \in \{(\f,\b), (\b,\f)\}$, a
  unary predicate $P$, a $C \in \Phi_0$ and terms $t_1$ and $t_2$, and
  we assume that $y$ occurs in $t_1$, and either $y$ does not occur in
  $t_2$ or $h_y(t_1) < h_y(t_2)$.  By the definition of $h_y(\phi)$
  and Axiom (F5), it suffices to prove $(\ast)$ with each of the
  atomic formulas $P (\g (t_1))$, $\g (t_1) = t_2$, $\g (t_1) <_C t_2$
  and $t_2 <_C \g (t_1)$ in place of $\phi$.  \medskip

  \noindent\textbf{Case 1:} $\phi$ is $P (\g (t_1))$.  By Axioms
  (F7)--(F9), the formula $\phi$ is equivalent in $T'$ to $\psi$,
  where $\psi$ is the formula depending on $P$ defined as follows:
  \begin{itemize}
  \item if $P \in \phiopen$ or $P$ is $E_F$ for some $F \in \phiopen$,
    then $\psi$ is
    \begin{equation*}
      \bigvee_{D \in \phitrans,\ P = D^\h} D (t_1) \vee \bigvee_{d \in
        \phisingle,\ P = d^\h} t_1 = d;
    \end{equation*}
  \item if $P \in \phitan$, then $\psi$ is the formula $t_1 = \h
    (e_P)$;
  \item if $P \in \phitrans$, then $\psi$ is the formula $E_{P^\h}
    (t_1)$.
  \end{itemize}
  In each case of $\psi$ above, we have $h_y(\psi) < h_y(\phi)$, as
  required. \medskip

  \noindent\textbf{Case 2:} $\phi$ is $\g (t_1) = t_2$.
  Then by Axioms (F5), (F7)--(F9) and (F13) the formula $\phi$ is
  equivalent in $T'$ to $\psi$, where $\psi$ is the conjuction of the
  formulas
  \begin{renumerate}
  \item ${\displaystyle t_2=s \ \vee \bigvee_{C \in \Phi_1} C(t_2) \
      \vee \bigvee_{c \in \phisingle} t_2=c \ \vee \bigvee_{C \in
        \phitan} t_2=e_C}$,
  \item $t_2 = c \to t_1 = \h (c)$ for each constant $c$ different from
    $s$,
  \item ${\displaystyle t_2=s \rightarrow \Bigg((t_1=s) \vee}$ \\
    ${\displaystyle \bigvee_{C \in \phiopen} \Big(E_C(t_1) \wedge
      \bigwedge_{D \in S_C}\neg(r^\h_D <_C t_1 <_C s^\h_D) \wedge
      \bigwedge_{c \in \phisingle}(\neg t_1=\h (c))\Big) \vee}$ \\
    ${\displaystyle \bigvee_{C \in \phitan} (\g (e_C) <_C t_1 \le_C e_C
      \vee e_C \le_C t_1 <_C \g (e_C)) \wedge \g (e_C)=s)\Bigg)}$ \\ with
    $S_C := \{D \in \phitrans:\ D^\h = C\}$,
  \item $C(t_2) \rightarrow t_1=\h (t_2)$ for $C \in \Phi_1$.
  \end{renumerate}
  If $y$ does not occur in $t_2$, then $h_y(\psi) < h_y(\phi)$; so we
  assume that $y$ occurs in $t_2$.  In this case, the only atomic
  subformula $\xi$ of $\psi$ with $h_y^1(\xi) = 1$ is $t_1 = \h (t_2)$,
  and $h_y(t_1 = \h (t_2)) = (1,h_y(t_1)) < (1, h_y(\g (t_1))) =
  h_y(\phi)$ by hypothesis, so $h_y(\psi) < h_y(\phi)$ as well.
  \medskip

  \noindent\textbf{Case 3:} $\phi$ is $\g (t_1) <_C t_2$.
  There are various subcases depending on $C$.
  \begin{itemize}
  \item If $C \in \phitrans$, we write $D:=C^\h$; then by Axioms (F8)
    and (F13) the formula $\phi$ is equivalent in $T'$ to $\psi$,
    where $\psi$ is the conjunction of the formulas
    \begin{equation*}
      (C(t_2) \vee t_2=\max(C)) \wedge ((E_D(t_1) \wedge r^\h_C <_D t_1 <_D
      r^\h_C) \vee t_1=\h(\min(C)))
    \end{equation*}
    and
    \begin{equation*}
      (E_D(t_1) \wedge r^\h_C <_D t_1 <_D r^\h_C) \rightarrow (t_1<_D \h
      (t_2) \vee t_2=\max(C)).
    \end{equation*}
  \item If $C \in \phiopen$, then by Axioms (F2), (F9), (F10), (F12)
    and (F13) the formula $\phi$ is equivalent in $T'$ to $\psi$,
    where $\psi$ is the conjunction of the formulas
    \begin{renumerate}
    \item ${\displaystyle \bigvee_{D \in \phitrans,\ D^\g=C}D(t_1) \vee
        \bigvee_{d \in \phisingle,\ P = d^\h}t_1=d}$,
    \item $\big(C(t_2) \wedge \neg E_C(t_2) \wedge
      E_C(\g (t_2))\big) \vee \big(C(t_2) \wedge \neg E_C(t_2) \wedge
      E_C(\h (t_2))\big) \vee E_C(t_2) \vee \big(t_2=\max(C)\big)$,
    \item $(D(t_1) \wedge E_C(t_2)) \rightarrow ((r^\g_D <_C t_2 <_C s^\g_D
      \wedge t_1<_D\h (t_2)) \vee (s^\g_D \leq_C t_2))$ for each $D
      \in \phitrans$ with $D^\g=C$,
    \item $(D(t_1) \wedge \neg E_C(t_2) \wedge E_C(\g (t_2))) \rightarrow
      ((r^\g_D <_C \g (t_2) <_C s^\g_D \wedge t_1 <_D\h (t_2)) \vee
      (s^\g_D \leq_C \g (t_2))$ for each $D \in \phitrans$ with
      $D^\g=C$,
    \item $(D(t_1) \wedge \neg E_C(t_2) \wedge E_C(\h (t_2))) \rightarrow
      ((r^\g_D <_C \h (t_2) <_C s^\g_D \wedge t_1 \leq_D \h(\h
      (t_2))) \vee (s^\g_D \leq_C \h (t_2)))$ for each $D \in \phitrans$
      with $D^\g=C$,
    \item $t_1=d \rightarrow \g d <_C t_2$ for $d \in \phisingle$ with
      $P = d^\h$.
    \end{renumerate}
  \item If $C \in \phitan$, then by Axioms (F2) and (F7) the formula
    $\phi$ is equivalent in $T'$ to $\psi'$, where $\psi'$ is
    \begin{equation*}
      (C (t_2) \vee t_2 = \max(C)) \wedge \big((t_1 = \h (e_C) \wedge e_C
      <_C t_2) \vee \g (t_1) = \min(C)\big).
    \end{equation*}
    In this case we let $\psi$ be the formula obtained from $\psi'$ by
    replacing the subformula $\g (t_1) = \min(C)$ by the corresponding
    formula obtained in Case 2.
  \end{itemize}
  We leave it to the reader to verify that $h_y(\psi) < h_y(\phi)$ in
  each of these subcases.
  \medskip

  \noindent\textbf{Case 4:} $\phi$ is $t_2 <_C \g (t_1)$.  This case is
  similar to Case 3; we leave the details to the reader.
\end{proof}

\begin{prop}
  \label{nice1} 
  Let $\phi(x,y)$ be a quantifier-free formula.  Then there is a
  minimal $y$-order formula $\psi(x,y)$ such that $\phi$ is
  equivalent in $T'$ to $\psi$.
\end{prop}

\begin{proof}
  By Corollary \ref{fsimp} and Lemma \ref{height_simplification}, we
  may assume that $\phi$ is a $y$-order formula such that $h^2_y(\phi)
  \le 1$.  By Lemma \ref{unmix}, there is a $y$-order formula
  $\psi'(x,y)$ such that $\phi$ is equivalent in $T'$ to $\psi'$,
  $\psi'$ contains no mixed terms and $h_y(\psi) \le h_y(\phi)$.

  In particular, for every binary atomic subformula $\eta$ of $\psi'$
  in which both terms contain $y$, one of the terms is $y$ itself and
  the other is either $\f^m(y)$ or $\b^m (y)$ for some $m = m(\eta) \in
  \NN$.  We now replace each such binary atomic subformula $\eta$ of
  $\psi'$ with $m(\eta) > 1$ by the formula $\eta'$ defined as
  follows:
  \begin{itemize}
  \item if $\eta$ is $y=\g^m (y)$ with $\g \in \{\f,\b\}$, then
    $\eta'$ is the disjunction of the formulas $y=c \wedge \g^m(c)=c$,
    for each constant symbol $c$, and $C(\g^m(y)) \wedge R_C(y,y)$,
    for each $C \in \Phi_1$;
  \item if $\eta$ is $y <_C \g^m (y)$ with $\g \in \{\f,\b\}$, then
    $\eta'$ is $B^\g_{m,C}(y,y)$;
  \item if $\eta$ is $\g^m (y) <_C y$ with $\g \in \{\f,\b\}$, then
    $\eta'$ is $S^\g_{m,C}(y,y)$.
  \end{itemize}
  We also replace each occurrence of $y=y$ by $s=s$ and each
  occurrence of $y <_C y$ by $s \ne s$, and we denote by $\psi''$ be
  the resulting formula.  Clearly $h_y(\psi'') \le h_y(\psi')$, and
  every binary atomic subformula of $\psi''$ in which both terms
  contain $y$ is of the form $G(y,y)$ for some $G \in \la(\Psi)
  \setminus \la(\Phi)$.  Moreover by Axioms (D1)--(D4), (D5)$_\nu$ and
  (D6)$_\nu$, the formula $\psi'$ is equivalent in $T'$ to $\psi''$.

  Next, we replace each subformula of $\psi''$ of the form $G(y,y)$,
  where $G \in \la(\Psi) \setminus \la(\Phi)$, by the corresponding
  minimal $y$-order formula $\psi(y)$ obtained in Lemma
  \ref{rsimp}(1).  If $\psi'''$ is the resulting $y$-order formula,
  then $\psi''$ is equivalent in $T(\Psi)$ to $\psi'''$ and
  $h_y^1(\psi''') = 0$.

  Finally by Lemmas \ref{height_simplification} and \ref{unmix}, there
  is a $y$-order formula $\psi$ such that $h_y(\psi) \le (0,1)$,
  $\psi$ contains no mixed terms and $\psi$ is equivalent in $T'$ to
  $\psi'''$.
\end{proof}

Finally, note that 
\begin{equation*}
  T(\Phi) \cup \{C(y)\} \models \neg E_C(y) \leftrightarrow \big(C (\f
  (y)) \vee C (\b (y))\big) 
\end{equation*}
for each $C \in \phiopen$, by Axioms (F5), (F10) and (F12).  Hence,
for each $C \in \phiopen$ and each $\g \in \{\f,\b\}$, we put
$T_{C,\g} := T(\Psi) \cup \{C(y) \wedge C (\g (y))\}$; by the previous
proposition, it remains to reduce quantifier-free formulas in each
$T_{C,\g}$.  It turns out, however, that we cannot entirely reduce to
minimal $y$-order formulas in these situations.

Instead, given $\g \in \{\f,\b\}$, we call a formula $\phi$
\textbf{$\g$-almost minimal} if $\phi$ is quantifier-free, the
only subterms of $\phi$ containing $z$ are $z$ and $\g (z)$ and every
binary atomic subformula $A(t_1,t_2)$ of $\phi$ is such that at most one
of $t_1$ and $t_2$ contains $z$.

\begin{prop}
  \label{nice2} 
  Let $\phi(x,y)$ be a quantifier-free formula, $C \in \phiopen$ and
  $\g \in \{\f,\b\}$.  Then there is a $\g$-almost minimal $y$-order
  formula $\psi_{C,\g}(x,y)$ such that $\phi$ is equivalent in
  $T_{C,\g}$ to $\psi_{C,\g}$.
\end{prop}

\begin{proof}
  By Corollary \ref{fsimp} and Lemma \ref{unmix}, we may assume that
  $\phi$ is a $y$-order formula containing no mixed terms.  On the
  other hand, we have $T \models \iota(\f (y))$ and $T \models
  \iota(\b (y))$ by Axiom (F5).  Let $\eta(x,y)$ be an atomic
  subformula of $\phi$; it suffices to show that there is a
  $\g$-almost minimal $y$-order formula $\xi_\eta(x,y)$ such that
  $\eta$ and $\xi_\eta$ are equivalent in $T_{C,\g}$.  If $h_y^2(\eta)
  = 0$, there is nothing to do, so we assume $h_y^2(\eta) > 0$, and we
  distinguish two cases to define $\xi_\eta$.  \medskip

  \noindent\textbf{Case 1:} $h_y^2(\eta) > 1$.  We first replace each
  occurrence of $\g (y)$ in $\eta$ by a new variable $z$ and each
  occurrence of $\h (y)$ in $\eta$ by $\h (z)$.  Denote the resulting
  atomic formula by $\eta'(x,z)$; by Axiom (F12), $\eta'(x,\g (y))$ is
  equivalent in $T_{C,\g}$ to $\eta(x,y)$.  By Proposition
  \ref{nice1}, the formula $\eta'(x,z)$ is equivalent in $T'$ to a
  minimal $z$-order formula $\eta''(x,z)$.  Since $T(\Psi) \models
  \iota(\g (y))$, it follows that $\eta$ is equivalent in $T_{C,\g}$ to
  the $\g$-almost minimal $y$-order formula $\xi_\eta$ given by
  $\eta''(x,\g (y))$.  \medskip

  \noindent\textbf{Case 2:} $h_y^2(\eta) = 1$.  In this case, we take
  $\xi_\eta$ equal to $\eta$ if $\eta$ contains a unary predicate
  symbol; so we assume that $\eta$ is a binary atomic formula
  $A(t_1,t_2)$.  If $\eta$ is $y = y$, we take $\xi_\eta$ to be $s=s$,
  and if $\eta$ is $y <_D y$ for some $D \in \Phi_0$, we take
  $\xi_\eta$ to be $s \ne s$; so we also assume from now on that
  $\max\{h_y^2(t_1), h_y^2(t_2)\} > 1$.  By Axiom (F5), the formulas
  $y = \g^m(y)$, $y = \h^m (y)$, $y <_D \g^m (y)$, $y <_D \h^m(y)$,
  $\g^m(y) <_D y$ and $\h^m(y) <_D y$, for $m > 0$ and $D \in \Phi_0
  \setminus \{C\}$, are all equivalent in $T_{C,\g}$ to $s \ne s$, so
  we are left with four subcases:
  \begin{renumerate}
  \item if $\eta$ is $y <_C \g^m (y)$ for some $m>0$, then we let
    $\eta'$ be the formula $(y <_C \g (y) \wedge C (\g^m (y)) \wedge
    R_C(\g (y)) \g y) \vee B^\g_{m-1,C}(\g (y), \g (y))$;
  \item if $\eta$ is $y <_C \h^m (y)$ for some $m>0$, then we let
    $\eta'$ be the formula $(y <_C \g (y) \wedge C (\h^m (y)) \wedge
    R_C (\g (y), \g (y))) \vee B^\h_{m,C}(\g (y), \g (y))$;
  \item if $\eta$ is $\g^m(y) <_C y$ for some $m>0$, then we let
    $\eta'$ be the formula $(\g (y) <_C y \wedge C (\g^m (y)) \wedge
    R_C (\g (y), \g (y))) \vee S^\g_{m-1,C} (\g (y), \g (y))$;
  \item if $\eta$ is $\h^m (y) <_C y$ for some $m>0$, then we let
    $\eta'$ be the formula $(\g (y) <_C y \wedge C (\h^m (y)) \wedge
    R_C (\g (y), \g (y))) \vee S^\h_{m,C} (\g (y), \g (y))$.
  \end{renumerate}
  We claim that $\eta$ and $\eta'$ are equivalent in $T_{C,\g}$.  We
  prove this for Case (i); the other cases are similar and left to the
  reader.  Let $b \in M$ be such that $\M \models C(b) \wedge C (\g
  (b))$.  Assume that $\M \models b <_C \g^m (b) \wedge \neg
  B^\g_{m-1,C}(\g (b), \g (b))$.  Then $\g^m (b) \in E_C$ and $\g^m
  (b) \le_C \g (b)$ by Axioms (F2) and (F5).  Hence $b <_C \g (b)$, so
  $\M \models \phi^\f(b, \g (b))$ by Axioms (F10) and (F12), which
  implies $\g^m (b) = \g (b)$ as required.  Conversely, assume first
  that $\M \models b <_C \g (b) \wedge C (\g^m (b)) \wedge R_C (\g
  (b), \g (b))$; then $b <_C \g^m (b)$ by Axioms (D2) and (F14).  Now
  assume that $\M \models B^\g_{m-1,C} (\g (b), \g (b))$; then $\g (b)
  <_C \g^m (b)$ by Axiom (D3), and hence $b <_C \g^m (b)$ by Axioms
  (F10) and (F12).

  Finally, by Proposition \ref{nice1}, the formulas $B^\g_{k,C}(z,z)$,
  $S^\g_{k,C}(z,z)$, $C (\g^{k} (z)) \wedge R_C(z,z)$ and $C (\h^{k}
  (z)) \wedge R_C(z,z)$ are each equivalent in $T'$ to minimal $z$-order
  formulas.  It follows from the claim that we are left to dealing
  with Subcases (i)--(iv) for $m=1$.  But by Axioms (F5), (F10) and
  (F12) we have $T_{C,\g} \models \neg C (\h (y))$.  Hence $T_{C,\g}
  \models \neg \phi^\h_C(y,\h (y))$, so from Axioms (F10) and (F12) we
  get $T_{C,\g} \models \phi^\g_C(y,\g y)$.  Therefore, $y <_C \g (y)$
  is equivalent in $T_{C,\g}$ to $s=s$ if $\g$ is $\f$, and to $\neg
  s=s$ if $\g$ is $\b$; the other subcases follow similarly.
\end{proof}

The previous two propositions allow us to reduce the problem of
eliminating quantifiers in $T(\Psi)$ to that of eliminating
quantifiers in two simpler theories:  for $C \in \Phi_1 \cup \phitan$
we let $\la_C$ be the language $\{<_C,\min(C),\max(C)\}$ and
$T_C$ be the $\la_C$-theory consisting of the universal
closures of
\begin{lenumerate}{A}{0}
\item the sentences stating that $<_C$ is a dense linear ordering on
  $C$, together with the formula $x = \min(C) \vee x = \max(C) \vee
  \min(C) <_C x <_C \max(C)$.
\end{lenumerate}
For $C \in \phiopen$ we let $\la_C$ be the language
$\{<_C,\pi_C,E_C,\min(C),\max(C)\}$, where $\pi_C$ a unary function
symbol, and we let $T_C$ be the $\la_C$-theory consisting of
the universal closures of (A1) as well as
\begin{lenumerate}{B}{0}
\item the formula $E_C (\pi_C (x)) \wedge (E_C(x) \to \pi_C (x) = x)$;
\item the formula $\pi_C (x) <_C x \to \neg\exists y (E_C(y) \wedge
  \pi_C (x) <_C y <_C x)$;
\item the formula $x <_C \pi_C (x) \to \neg\exists y (E_C(y) \wedge x
  <_C y <_C \pi_C (x))$;
\item the sentences stating that for every $x \in E_C$, the
  restriction of $<_C$ to the set $\set{y:\ \pi_C (y) = x}$ is a dense
  linear ordering without endpoints.
\end{lenumerate}

A routine application of a quantifier elimination test such as Theorem
3.1.4 in \cite{mark2} gives the following result; we leave the
details to the reader.

\begin{prop}
  \label{lexQE}
  For each unary predicate symbol $C$ of $\la(\Phi)$, the theory
  $T_C$ admits quantifier elimination in the language $\la_C$. \qed
\end{prop}

\begin{thm}
  \label{main} 
  The theory $T(\Psi)$ admits quantifier elimination.
\end{thm}

\begin{proof}
  Let $\phi(x,y)$ be a quantifier-free formula; we show that $\exists
  y \phi(x,y)$ is equivalent in $T(\Psi)$ to a quantifier-free
  formula.  First, note that $\exists y \phi(x,y)$ is equivalent
  in $T(\Psi)$ to the disjunction of the formulas
  \begin{enumerate}
  \item $\phi(x,c)$ for each constant $c$;
  \item $\exists y (C(y) \wedge \phi(x,y))$ for each $C \in \Phi_1 \cup
    \phitan$;
  \item $\exists y (C(y) \wedge C \g (y) \wedge \phi(x,y))$ for each $C
    \in \phiopen$ and each $\g \in \{\f,\b\}$.
  \end{enumerate}
  We deal with each disjunct separately; since formulas of type (1)
  are trivial to handle, we deal with types (2) and (3).  \medskip

  \noindent\textbf{Type (2):}  Let $C \in \Phi_1 \cup \phitan$.  Since
  $T(\Psi) \models C(y) \to \iota(y)$, we may assume by Proposition
  \ref{nice1} that $\phi$ is a minimal $y$-order formula.  Without
  loss of generality, we may also assume that $\phi$ is a conjunction
  of atomic formulas, that $y$ occurs in each of the atomic
  subformulas of $\phi$ and, by Axiom (F1), that $\phi$ contains only
  the relation symbols $=$ and $<_C$.  Let $t_1, \dots, t_k$ be all
  maximal subterms of $\phi$ that do not contain $y$, and let
  $\phi'(z_1, \dots, z_k,y)$ be the formula obtained from $\phi$ by
  replacing each $t_i$ by a new variable $z_i$.  Then $\phi'$ is a
  $<_C$-formula without parameters; by Proposition \ref{lexQE}, there
  is a quantifier-free $\la_C$-formula $\psi'(z_1, \dots, z_k)$
  such that $\exists y \phi'$ and $\psi'$ are equivalent in
  $T_C$.  Let $\psi(x)$ be the $\la(\Psi)$-formula obtained from
  $\psi'$ by replacing each $z_i$ by $t_i$; then $\exists y \phi$ and
  $\psi$ are equivalent in $T(\Psi)$, as required.  \medskip

  \noindent\textbf{Type (3):}  Let $C \in \phiopen$ and
  $\g \in \{\f,\b\}$; by Proposition \ref{nice2}, we may assume that
  $\phi$ is a $\g$-almost minimal $y$-order formula.  Without loss of
  generality, we may also assume that $\phi$ is a conjunction of
  atomic formulas, that $y$ occurs in each of the atomic subformulas
  of $\phi$ and, by Axiom (F1), that $\phi$ contains only the relation
  symbols $=$, $<_C$ and $E_C$.  Let $t_1, \dots, t_k$ be all maximal
  subterms of $\phi$ that do not contain $y$, and let $\phi'(z_1,
  \dots, z_k,y)$ be the formula obtained from $\phi$ by replacing each
  $t_i$ by a new variable $z_i$.  Note that $\phi'$ contains no
  parameters.  Arguing as for Type (2), it now suffices to find a
  quantifier-free formula $\psi'(z_1, \dots, z_k)$ equivalent in
  $T(\Psi)$ to $\exists y \phi'(z_1, \dots, z_k,y)$.

  To do so, we let $\pi_C$ be a new unary function symbol and let
  $T(\Psi)_C$ be the theory $T(\Psi)$ together with the universal closure
  of the formula
  \begin{multline*}
    y = \pi_C (x) \leftrightarrow \Big(\big(E_C(x) \wedge y = x\big)
    \\ \vee \big(C(x) \wedge C (\f (x)) \wedge y = \f (x)\big) \vee
    \big(C(x) \wedge C (\b (x)) \wedge y = \b (x)\big)\Big).
  \end{multline*}
  Since $T(\Psi)_C$ is an extension by definitions of $T(\Psi)$ in the
  sense of \cite[Section 4.6]{shoenfield}, it suffices to find a
  quantifier-free $\la(\Psi)$-formula $\psi'(z_1, \dots, z_k)$
  equivalent in $T(\Psi)_C$ to $\exists y \phi'(z_1, \dots, z_k,y)$.

  Let $\phi''$ be the $\la_C$-formula obtained from $\phi'$ by
  replacing each occurrence of $\g (y)$ by $\pi (y)$; then $\phi'$ and
  $\phi''$ are equivalent in $T(\Psi)_C$.  Since $T(\Psi)_C \models T_C$,
  there is by Proposition \ref{lexQE} a quantifier-free
  $\la_C$-formula $\psi''(z_1, \dots, z_k)$ that is equivalent in
  $T(\Psi)_C$ to $\exists y \phi''(z_1, \dots, z_k,y)$; without loss of
  generality, we may assume that the only subterms of $\psi''$ are
  $z_i$ and $\pi z_i$ for $i = 1, \dots, k$.

  Finally, we let $\psi'$ be the $\la(\Psi)$-formula obtained from
  $\psi''$ by replacing each atomic subformula $\eta$ of $\psi''$ by
  an $\la(\Psi)$-formula $\eta'$ determined as follows:
  \begin{renumerate}
  \item if $\eta$ is $E_C(\pi_C (z_i))$, we let $\eta'$ be $C (z_i)
    \wedge (E_C (z_i) \vee C (\f (z_i)) \vee C (\b (z_i)))$;
  \item if $\eta$ is $\pi_C (z_i) \ast z_j$ with $\ast \in \{=, <_C,
    >_C\}$, we let $\eta'$ be
    \begin{equation*}
      C (z_i) \wedge C (z_j) \wedge \left(\bigvee_{\g \in \{\f^0, \f,
          \b\}} E_C (\g (z_i)) \wedge \g (z_i) \ast z_j\right);
    \end{equation*}
  \item if $\eta$ is $\pi_C (z_i) <_C \pi_C (z_j)$ and $\ast \in \{=,
    <_C\}$, we let $\eta'$ be
    \begin{equation*}
      C (z_i) \wedge C (z_j) \wedge \left(\bigvee_{\g,\h \in \{\f^0, \f,
          \b\}} E_C (\g (z_i)) \wedge E_C (\h (z_j)) \wedge \g (z_i)
        \ast \h (z_j)\right); 
    \end{equation*}
  \end{renumerate}
  and if $\eta$ is not of one of the forms (i)--(iii) above, we let
  $\eta'$ be $\eta$.  This $\psi'$ is equivalent in $T(\Psi)_C$ to $\psi''$
  and is of the required form.
\end{proof}

\section{Consequences for the Model Theory of $T(\Psi)$}  
\label{consequences}

The quantifier elimination result established in the previous section
allows us to show that the theory $T(\Psi)$ is very well-behaved: it
is a theory of finite rank in the sense developed by Onshuus
\cite{alf1}.  

We first rephrase the results from the previous section.  For a flow
configuration $\Phi$, $C \in \phiopen$, $\M \models T(\Psi)$ and $x
\in E_C^\M$, we put
\begin{equation*}
  C_x^\M := \set{y \in C^\M:\ y=x \vee \f (y) = x \vee \b (y) = x}
\end{equation*}
and $\bar{C}^\M_x := C^\M_x \cup \{\f(x),\g(x)\}$.  The following
corollary implies Theorem C:

\begin{cor}
  \label{induced}
  Let $\Psi$ be a Dulac flow configuration and $\M \models T(\Psi)$.
  \begin{enumerate}
  \item For $C \in \Phi_1 \cup \phitan$, every definable subset
    of $C^\M$ is a finite union of points and open $<_C$-intervals
    with endpoints in $\bar{C}$.
  \item For $C \in \phiopen$ and $x \in E_C^\M$, every definable
    subset of $C_x^\M$ is a finite union of points and open
    $<_C$-intervals with endpoints in $\bar{C}^\M_x$.
  \end{enumerate}
\end{cor}

\begin{proof}
  This follows immediately from Theorem \ref{main}, Propositions
  \ref{nice1} and \ref{nice2} and Axioms (F2) and (F11).
\end{proof}

Below we use the terminology of rosy theories.

\begin{thm} 
  \label{rank2}
  Let $\Psi$ be a Dulac flow configuration and $T$ be any completion
  of $T(\Psi)$.  Then $T$ is rosy with $\urank(T) \le 2$.
\end{thm}

\begin{proof}
  Let $p(x)$ be a complete $1$-type in $T$, $\M \models T$ and $a \in
  M$ such that $\M \models p(a)$.  If $C(x) \in p$ for some $C \in
  \phitan \cup \Phi_1$, then by Proposition \ref{induced}(1) the type
  $p$ is determined by the $<_C$-order type of $x$ over the constants;
  hence $\urank(p) \le 1$.  If $C(x) \wedge \neg E_C(x) \in p$ for
  some $C \in \phiopen$, then by Proposition \ref{induced}(2) the type
  $p$ is determined by the $<_C$-order type $o(x)$ of $a$ over the
  constants and $\pi_C(a)$, where $\pi_C:C \into E_C$ is given by
  \begin{equation*}
    \pi_C(z) := 
    \begin{cases}
      z &\text{if } z \in E_C^\M, \\ \f(z) &\text{if } \f(z) \in
      E_C^\M, \\ \b(z) &\text{if } \b(z) \in E_C^\M.
    \end{cases}
  \end{equation*}
  Again by Proposition \ref{induced}(1), the type of $\pi_C(a)$ over the
  constants is determined by the $<_C$-order type of $\pi_C(a)$ over
  the constants. 

  Since $p$ either contains one of the above formulas or a formula $x
  = c$ for some constant symbol $c$, it follows from the Fact in the
  introduction that $\urank(T) \leq 2$.
\end{proof}

In fact, the $\urank$-rank in the previous theorem is actually equal
to 2:

\begin{prop}
  \label{minrank}
  Let $\Phi$ be a flow configuration and $\M \models T(\Phi)$, and
  assume that $\phiopen \ne \emptyset$.  Then $\urank(\M) \ge 2$.
\end{prop}

\begin{proof}
  Let $C \in \phiopen$.  Then by the example in the introduction, the
  theory of $(C,<_C,E_C)$ has $\urank$-rank at least two.  Hence
  $\urank(\M) \ge 2$.
\end{proof}

There is a certain converse to Theorem \ref{rank2} based on Remark
\ref{predulac}: we let $\Phi$ be a flow configuration and consider the
theory $T(\Phi)^+$ obtained by adding the universal closures of the
following formulas to $T(\Phi)'$ for each $C \in \phitrans$:
\begin{equation}  \label{DC}
  \begin{split}
    C(x) \to \exists y \left(\bar{C} (y) \wedge y =
      \inf\{z:\ x <_C z \wedge \fixbd_C(z)\}\right) \\
    C(x) \to \exists y \left(\bar{C} (y) \wedge y = \sup\{z:\ z <_C x
      \wedge \fixbd_C(z)\}\right).
  \end{split}
\end{equation}

\begin{expls}
  \label{poincare}
  (1)  Let $\Psi$ be a Dulac flow configuration.  Then any model $\M$
  of $T(\Psi)$ satisfies \eqref{DC}.

  (2) Let $\xi$ be a definable vector field on $\RR^2$, and let
  $\M_\xi$ be an $\la(\Phi_\xi)$-structure associated to
  $\xi$ as in Example \ref{omega_structure}.  Then $\M_\xi$
  satisfies \eqref{DC} by Corollary \ref{def_complete}, and by
  Remark \ref{predulac} the structure $\M_\xi$ can be expanded to a
  model $\M^+_\xi$ of $T(\Phi_\xi)^+$.
\end{expls}

Below, for each $\nu \in \NN$ we abbreviate the formula stating that
$\fixbd_C(x)$ defines a set with at most $\nu$ elements by
``$|\fixbd_C(x)| \le \nu$''.

\begin{prop}
  \label{converse}
  Let $\Phi$ be a flow configuration and $T$ be a completion of
  $T(\Phi)^+$, and assume that $\urank(T) \le 2$.  Then there is a
  $\nu \in \NN$ such that
  \begin{enumerate}
  \item $T \models |\fixbd_C(x)| \le \nu$;
  \item every model $\M$ of $\,T$ can be expanded to a model of
    $T(\Phi,\nu)$.
  \end{enumerate}
\end{prop}

\begin{proof}
  (1) Assume that $T \not\models |\fixbd_C(x)| \le \nu$ for any $\nu
  \in \NN$.  Then by model theoretic compactness, there are an $\M
  \models T$ and a $C \in \Phi_1$ such that the set $\fixbd_C(M)$ is
  infinite; we may assume that $\M$ is $\aleph_1$-saturated.  Moreover
  by Axiom (F8), we may assume that $C \in \phitrans$.  Also, by Axiom
  (F8) and an argument as in the proof of Proposition \ref{minrank},
  it suffices to find a $d \in C^{\M}$ such that $\urank(d) \ge 2$.

  Since $\M$ is $\aleph_1$-saturated, there is an interval $I
  \subseteq C^\M$ such that $I \cap \acl(\emptyset) = \emptyset$ and
  $I \cap \fixbd_C(M)$ is infinite.  By \eqref{DC} and since
  $\fixbd_C(M)$ is nowhere dense, there is a $c \in I \setminus
  \fixbd_C(M)$ such that the elements $a:= \sup\set{x \in I:\ x <_C c
    \wedge \fixbd_C(x)}$ and $b:= \inf\set{x \in C:\ a <_C x \wedge
    \fixbd_C(x)}$ exist in $I$.  Then $a <_C b$, $a,b \notin
  \acl(\emptyset)$, $b \in \dcl(a)$ and
  \begin{equation*}
    \M \models a <_C b \wedge \fixbd_C(a) \wedge \neg\exists x (C(x)
    \wedge a <_C x  <_C b \wedge \fixbd_C(x)). 
  \end{equation*}
  It follows that the formula $\phi(x):=a<_Cx<_Cb$ strongly divides
  over $\emptyset$; hence $\urank(d) \geq 2$ for some $d \in C^\M$, as
  required.

  Part (2) follows from Proposition \ref{dulacprop} and part (1).
\end{proof}

We can now prove our restatement of Dulac's Problem:

\begin{proof}[Proof of Theorem B]
  (1) If $\xi$ has finitely many boundary cycles, then by
  Pro\-position \ref{dulacprop} the structure $\M_\xi$ can be
  expanded into a model $\M^D_\xi$ of $T(\Phi_\xi,\nu)$ for some
  $\nu \in \NN$.  Since $(\Phi_\xi)_{\text{open}} \ne \emptyset$,
  it follows that $2 \le \urank(\M_\xi) \le \urank(\M^D_\xi) \le
  2$ by Proposition \ref{minrank} and Theorem \ref{rank2}.
  
  Conversely, if $\urank(\M_\xi) = 2$ then by Proposition
  \ref{converse}, the structure $\M_\xi$ can be expanded into a
  model of $T(\Phi_\xi,\nu)$ for some $\nu \in \NN$, so by Example
  \ref{limit_cycle_set} the vector field $\xi$ has finitely many
  boundary cycles.

  Part (2) follows from part (1) and Poincar\'e's Theorem
  \cite{poincare} (see also \cite[p. 217]{perko}).  The ``moreover''
  clause follows from part(1) and Theorem \ref{rank2}.
\end{proof}

\section{Final questions and remarks}  \label{final}

\begin{enumerate}
\item In the situation of Theorem B, is it possible for
  $\M_\xi$ to be rosy of $\urank$-rank strictly greater than $2$?
\item Can a restatement of Hilbert's 16th Problem be obtained in the
  spirit of Theorem B?  

  A na\"ive approach to this question is as follows: Let $\{\xi_a:\
  a \in A\}$ be a family of vector fields on $\RR^2$ definable in $\R$.
  Since the arguments in Sections \ref{rolle} through
  \ref{progression_map} are uniform in parameters, we may assume that
  there is a flow configuration $\Phi$ such that $\Phi_{\xi_a} =
  \Phi$ for all $a \in A$.  In this situation, one can readily
  reformulate the theory $T(\Phi)$ for the parametric situation; and
  if one also assumes the existence of a uniform bound $\nu \in \NN$
  on the number of boundary cycles of each $\xi_a$, such a
  reformulation extends to $T(\Phi,\nu)$.  We suspect that under the
  latter assumption, the corresponding theory is rosy of $\urank$-rank
  $3$; however, this does not appear to us to be a completely trivial
  generalization of the results in Section \ref{consequences}, and we
  plan to pursue it in a future project.
\item The structure $\M_\xi^D$ in Example
  \ref{dulac_omega_structure} does not define any algebraic operations
  (by Theorem \ref{main}).  Assume here that $S(\xi) = \emptyset$;
  is it possible to expand $\M_\xi^D$ by some (or all) of the sets
  definable in the original o-minimal structure $\R$ without
  increasing the $\urank$-rank?  We know very little about this
  question.  However, if (a) the $x$-axis, the projection from $\RR^2$
  onto the $x$-axis, and both addition and multiplication are
  definable in an expansion $\M'$ of $\M_\xi^D$, and if (b) the
  expansion $\M'$ still has $\urank$-rank two, then $\M'$ (and hence
  $\M_\xi^D$) would be definable in an o-minimal structure.  (The
  assumption that $\M'$ has $\urank$-rank two is necessary here.)
  Thus, question (3) is related to the following question:
\item Is the structure $\M_\xi^D$ of Example
  \ref{dulac_omega_structure} definable in some o-minimal expansion of
  the real field?  
\item Consider a Dulac flow configuration $\Psi$ and $\M \models
  T(\Psi)$.  Corollary \ref{induced}, Theorem \ref{rank2} and their
  respective proofs may be loosely interpreted as indicating that $\M$
  is built-up from sets $D \subseteq M$ on which the induced structure
  is o-minimal.  Is there a theory of structures built-up from sets
  with induced o-minimal structure, say in the spirit of Zilber's
  results on the fine structure of uncountably categorical theories
  \cite{zilber}?  
\end{enumerate}

\end{document}